

\magnification=\magstephalf
\hoffset=0.5 truein             

\input amstex
\loadeusm
\loadbold

\documentstyle{amsppt}
\refstyle{A}

\define\cE{{\Cal E}}\define\cF{{\Cal F}}
\define\cI{{\Cal I}}
\define\cJ{{\Cal J}}

\define\bbR{{\Bbb R}}\define\bbC{{\Bbb C}}
\define\bbP{{\Bbb P}}
\define\w{{\mathchoice{\,{\scriptstyle\wedge}\,}{{\scriptstyle\wedge}}
{{\scriptscriptstyle\wedge}}{{\scriptscriptstyle\wedge}}}}
\define\lefthook{\mathbin{
\hbox{\vrule height1.4pt width4pt depth-1pt
\vrule height4pt width0.4pt depth-1pt}}}

\define\abold{{\bold a}}\define\bbold{{\bold b}}
\define\ebold{{\bold e}}\define\sbold{{\bold s}}
\define\rbold{{\bold r}}\define\ubold{{\bold u}}\define\pbold{{\bold p}}
\define\Pbold{{\bold P}}\define\Qbold{{\bold Q}}\define\Rbold{{\bold R}}
\define\Xbold{{\bold X}}\define\Hbold{{\bold H}}
\define\euso{{\frak so}}
\define\SO{\text{SO}}\define\Or{\text{O}}
\define\ts{\textstyle }

\define\>{\medspace}
\define\Lie{{\eusm L}}

\define\alphabold{{\boldsymbol\alpha}}
\define\rhobold{{\boldsymbol\rho}}
\define\etabold{{\boldsymbol\eta}}
\define\sigmabold{{\boldsymbol\sigma}}
\define\thetabold{{\boldsymbol\theta}} 
\define\Thetabold{{\boldsymbol\Theta}}

\topmatter
\title 
   Harmonic morphisms with fibers of dimension one
\endtitle
\author
   Robert L. Bryant
\endauthor
\affil
   Duke University
\endaffil
\address
   Department of Mathematics,
   Duke University,
   PO Box 90320,
   Durham, NC 27708-0320
\endaddress
\date
   January 13, 1997
\enddate

\email
   bryant\@math.duke.edu
\endemail

\keywords
harmonic maps, harmonic morphisms, exterior differential systems
\endkeywords

\subjclass
58E20,  
58A15,  
53C12   
\endsubjclass

\thanks
This research was begun during a visit to IMPA in Rio de Janeiro
in July 1996 and was inspired by questions raised during the 
International Conference on Differential Geometry held at
IMPA during that month.  The article was written during
a visit to the Institute for Advanced Study in Princeton. 
The author would like to thank IMPA and the IAS for their hospitality 
and also to acknowledge support from the National Science Foundation 
through grant DMS--9505125.
\endgraf
The current version of this preprint is at
{\tt file://www.duke.edu/\lower.5ex\hbox{\tt\char'176}bryant/HarMorphs.dvi}\ .  
\endthanks

\abstract
The harmonic morphisms~$\phi:M^{n+1}\to N^n$ are studied using
the methods of the moving frame and exterior differential
systems and three main results are achieved.  

The first result is a local structure theorem for such maps in the case 
that~$\phi$ is a submersion, in particular, a normal form is found for 
all such~$\phi$ once the metric on the target manifold~$N$ is specified.  

The second result is a finiteness theorem, 
which says, in a certain sense, that, when $n\ge3$, the set of 
harmonic morphisms with a given Riemannian domain~$\bigl(M^{n+1},g\bigr)$
is a finite dimensional space.  

The third result is the explicit classification when~$n\ge3$ 
of all local and global harmonic morphisms with domain~
$\bigl(M^{n+1},g\bigr)$, a space of constant curvature.  
\endabstract

\endtopmatter

\document

\head 0. Introduction \endhead

A smooth map~$\phi:M\to N$ between Riemannian manifolds is said to be
a {\it harmonic morphism\/} if, for any harmonic function~$f$ on any open 
set~$V\subset N$, the pullback~$f\circ\phi$ is a harmonic function 
on~$\phi^{-1}(V)\subset M$. 

By a simple argument (see~\S1), any non-constant harmonic morphism~$\phi:M\to 
N$ between connected Riemannian manifolds must be a submersion away from a set 
of measure zero in~$M$.  Thus, a necessary condition for the existence of a 
non-constant harmonic map~$\phi:M\to N$ is that~$\dim M\ge\dim N$. 

When the dimension of~$N$ is~$1$, so that $N$ can be regarded, at least 
locally, as~$\bbR$ with its standard metric, a map~$\phi:M\to N$ is a 
harmonic morphism if and only if it is a harmonic function in the usual sense.  
Thus, at least locally, there are many harmonic morphisms from~$M$ to~$N^1$.  

However, when the dimension of~$N$ is greater than~$1$, the condition 
of being a harmonic morphism turns out to be much more restrictive, 
being essentially equivalent to an overdetermined system of~{\smc pde} for 
the map~$\phi$.  Thus, for generic Riemannian metrics on~$M$ 
and~$N$, one does not expect there to be any harmonic morphisms, 
even locally.  Moreover, in the case that there do exist harmonic
morphisms~$\phi:M\to N$ for given~$M$ and~$N$, one expects the analysis
of the overdetermined system that describes them to involve integrability
conditions and other features of overdetermined systems.

When both~$M$ and~$N$ have dimension~$2$, a harmonic morphism is simply a 
branched conformal mapping between Riemann surfaces and these are studied 
by classical methods of complex analysis and Riemann surface theory.  

When both manifolds have the same dimension~$n>2$, a non-constant harmonic 
morphism is a local homothety, i.e., up to a constant scale factor, $\phi$ 
is a local isometry.  

Thus, the interesting cases are when $\dim M>\dim N\ge2$. This article 
concerns the case when~$\dim M = \dim N +1$, i.e., when the dimension of the 
generic fiber of~$\phi$ is~1.  It contains three main results.

The first, Theorem~1, is a local structure theorem for harmonic morphisms 
whose fibers are curves.  This result
describes the possible Riemannian metrics~$g$ that can be defined on the 
domain~$M$ of a smooth mapping~$\phi:M\to N$ where $N$ is a smooth manifold
endowed with a fixed Riemannian metric~$h$ so that~$\phi:(M,g)\to(N,h)$
will be a harmonic morphism.  

The second result, Theorem~2,
is a general finiteness theorem for harmonic morphisms of corank one with
a given Riemannian domain~$(M^{n+1},g)$ where $n\ge 3$.  This result shows
that the set of such harmonic morphisms is, in a certain sense, finite
dimensional.  This result is in marked contrast to the case~$n=2$,
which has already been analyzed by Baird and Wood with the result that the 
locally defined harmonic morphisms with a given Riemannian domain~$(M^3,g)$ 
of constant sectional curvature depend on arbitrary functions (in the
sense of exterior differential systems).

The third result, Theorem~3, is a classification of the harmonic morphisms 
of corank one whose domain~$\bigl(M^{n+1},g\bigr)$ is a simply-connected,
complete Riemannian manifold of constant curvature and dimension $n{+}1\ge 4$.
It will be shown that there are exactly two types of such harmonic morphisms.  

The first type
can be thought of as a sort of metric quotient and is described as follows:
Let~$X$ be a Killing vector field on~$M$ with zero locus~$Z\subset M$ and 
suppose that the space~$N$ of integral curves of~$X$ in~$M\setminus Z$ can 
be given the structure of a smooth $n$@-manifold in such a way that the 
quotient map~$\phi:M\setminus Z\to N$ is a smooth submersion.
Then there exists a metric~$h$ on~$N$, unique up to a constant scale factor,
so that~$\phi:(M\setminus Z,g)\to(N,h)$ is a harmonic morphism.  (Sometimes
this map can be extended across~$Z$ as well after suitably extending~$N$,
see \S3.3.)
The second type is described as follows:  Let $N\subset M$ be a totally 
umbilic hypersurface, endowed with a constant multiple of the induced metric, 
denoted~$h$.  Let~$P\subset M$ be the focal set of~$N$, which consists of at 
most two points.  There is a canonical retraction~$\phi:M\setminus P\to N$ 
that retracts~$M\setminus P$ back to~$N$ along the geodesics normal to~$N$. 
Then~$\phi$ is a harmonic morphism.  The examples of this kind had already
appeared in the work of Gudmundsson~\cite{Gu1}.

The methods used are those of exterior differential systems and the
moving frame, both of which are well-adapted to the study of overdetermined
systems of {\smc pde}.

\remark{Acknowledgments}
It is a pleasure to thank John Wood, whose questions inspired this article 
and whose comments and and guide to the literature on harmonic morphisms 
were invaluable.
\endremark

\head 1. Harmonic morphisms via moving frames \endhead

This section is a self-contained treatment by moving frame calculations of
the the basic structure theory of harmonic morphisms.  It is intended to be 
readable by those familiar with either moving frame calculations or 
the fundamentals of harmonic morphisms. 

Its main purpose is to fix notation and to serve as a reference for the 
proofs in the later sections, which employ the moving frame.  Such a reference
is probably needed, as it appears that most of the current workers 
on harmonic morphisms do not use moving frames and so the translation of 
known results into this language may be helpful.  For more background
on the method of the moving frame, see~\cite{Sp}.  Most of the results
about harmonic maps and morphisms to be derived in this section  
can be found in the standard references on the subject, 
such as~\cite{EL1},~\cite{EL2}, or~\cite{W2}.

\subhead 1.1. Moving frame computations for harmonic morphisms \endsubhead
Let $M$ and $N$ be Riemannian manifolds of dimensions~$m$ and~$n$, 
respectively.  For simplicity, I assume that both~$M$ and $N$ are connected
throughout this article.   The summation convention will be used extensively, 
with the understood ranges
$$
\alignedat2 
1&\le a,b,c&&\le m,\\ 
1&\le i,j,k&&\le n.\\
\endalignedat
$$

\subsubhead 1.1.1. 
Coframings, connection forms, and structure equations
\endsubsubhead
Let $g$ be the metric on~$M$ and $h$ be the metric on~$N$.  
Let $U\subset M$ and $V\subset N$ be open sets with trivial tangent bundles.  
Then there exist smooth 
coframings~$\omega = \bigl(\omega_1,\dots,\omega_m\bigr)$ 
and~$\eta = \bigl(\eta_1,\dots,\eta_n\bigr)$ on~$U$ and~$V$ respectively, 
so that
$$
\aligned
g_{|U} & = \omega_1^2+\dots+\omega_m^2 = {\omega_a}^2\,\\
h_{|V} & = \eta_1^2+\dots+\eta_n^2 = {\eta_i}^2\,.\\
\endaligned
\tag1
$$

Corresponding to the chosen coframings on the respective
open sets, there exist unique $1$-forms~$\omega_{ab}=-\omega_{ba}$ 
and $\eta_{ij}=-\eta_{ji}$ that represent the Levi-Civita connections 
of the respective metrics and that are characterized by the {\it structure 
equations\/}
$$
d\omega_a = -\omega_{ab}\w\omega_b\,,
\qquad\qquad
d\eta_i = -\eta_{ij}\w\eta_j\,.
\tag2
$$

\subsubhead 1.1.2. 
Mappings and pullbacks
\endsubsubhead
Now suppose that $\phi:M\to N$ is a smooth map and that $U$ and
$V$ have been chosen so that $U\subset f^{-1}(V)$.  Then there
exist unique functions~$f_{ia}$ on~$U$ so that
$$
\phi^*(\eta_i) = f_{ia}\,\omega_a\,.
\tag3
$$

Because the chosen coframings are orthonormal, the {\it energy 
density\/} of the map~$\phi$ on~$U$ is given by
$$
\align
E(\phi)_{|U} &= f_{ia}f_{ia}\,|\omega_1\w\dots\w\omega_n|\\
             &={||\phi'||^2}_{|U}\,|\omega_1\w\dots\w\omega_n|\,.\\
\endalign
$$
This density is globally defined, independent of the local choice of~
$\omega$ or~$\eta$. When~$M$ is compact, integration of this density yields 
a functional called the {\it energy\/}~$\cE:C^\infty(M,N)\to\bbR$, 
namely
$$
\cE(\phi) = \frac12\int_M E(\phi).
$$

Adopt the convention that, for any differential form~$\psi$ on~$V$,
its $\phi$-pullback~$\phi^*(\psi)$ on~$U$ is denoted by
an overbar, i.e., $\overline\psi = \phi^*(\psi)$.  Since the map~$\phi$ 
will be fixed in this discussion, this should cause no
confusion.  Thus, (3) becomes 
$$
\overline{\eta_i} = f_{ia}\,\omega_a\,.
\tag$3'$
$$
(The reader of other sources on moving frame calculations should be 
aware that many authors simply drop the pullback notation entirely, writing
(3) in the even simpler form~$\eta_i = f_{ia}\,\omega_a$.  This has
caused considerable confusion in some cases, a confusion I hope to 
avoid.)

Taking the exterior derivative of (3) and using the structure 
equations~(2) yields
$$
\bigl(df_{ia}-f_{ib}\,\omega_{ba}+f_{ja}\,\overline{\eta_{ij}}\bigr)
\w\omega_a = 0.
$$
By Cartan's Lemma, there exist unique functions~$f_{iab}=f_{iba}$
on~$U$ so that
$$
df_{ia}=f_{ib}\,\omega_{ba}-f_{ja}\,\overline{\eta_{ij}}+f_{iab}\,\omega_b\,.
\tag4
$$
The {\it tension field\/} of~$\phi$ on~$U$ is the tensor field
$$
\tau(\phi) = f_{iaa}\>\ebold_i{\circ}\phi
$$
(where $(\ebold_1,\ldots,\ebold_n)$ is the dual orthonormal frame field 
in~$V$ to $\eta$) and is a section of the $\phi$-pullback to~$U$ of the 
tangent bundle of~$N$.  This tensor does not depend on the 
choice of~$\omega$ or~$\eta$ and so is globally defined on~$M$.

The map~$\phi$ is said to be {\it harmonic\/} if it satisfies the 
Euler-Lagrange equations for the energy functional~$\cE$. A calculation 
(see~\cite{EL2,(2.30)}) shows that~$\phi$ is harmonic if and only if its 
tension field vanishes.  A function~$f:N\to\bbR$ is said to be harmonic 
if it is harmonic as a map to~$\bbR$ endowed with its standard metric.

\subsubhead 1.1.3. 
Harmonic morphisms
\endsubsubhead
It is not generally true that the composition of harmonic maps is
harmonic, so there is no useful category whose objects are the Riemannian 
manifolds and whose morphisms are harmonic maps.  Nevertheless, some
harmonic maps turn out to have useful composition properties and these
are designated as harmonic morphisms.  In this subsection, this notion
will be explained and the equations for it derived.

Given a smooth function~$v$ on~$V$, there exist functions~$v_i$ 
so that~$dv = v_i\,\eta_i$.  Differentiating this relation gives 
$0 =\bigl(dv_i - v_k\,\eta_{ki}\bigr)\w\eta_i$.  By Cartan's Lemma,
there exist functions~$v_{ij}=v_{ji}$ so that  
$dv_i=v_k\,\eta_{ki}+v_{ij}\,\eta_j$.  The trace~$\Delta v =
-v_{ii}$ is the local expression with respect to the coframing~$\eta$ for 
the globally defined $h$@-Laplacian of~$v$, which is also the
tension field of~$v$ regarded as a mapping~$v:V\to\bbR$. 
Thus, $v$ is harmonic in~$V$ if and only if ~$\Delta v$ vanishes 
identically in~$V$.

Locally, there are many harmonic functions.
It is known~\cite{EL1} that, for any~$q\in V$, and any collection
of real numbers~$r_i$ and $r_{ij}=r_{ji}$ with $r_{ii}=0$, there is an open 
$q$@-neighborhood~$V'\subset V$ and a harmonic function~$v$ on~$V'$ so
that~$v_i(q)=r_i$ and $v_{ij}(q)=r_{ij}$.  This local `flexibility' will
be important below.

Roughly speaking, a map~$\phi:M\to N$ will be a harmonic `morphism' if 
$\phi$@-pullback carries harmonic functions on~$N$ to harmonic functions
on~$M$.  The precise definition is as follows:

\proclaim{Definition 1}
A map~$\phi:M\to N$ is a {\it harmonic morphism\/} if, for any 
harmonic function~$v$ on an open set~$V\subset N$, the function~$\overline{v} 
= \phi^*(v)$ is harmonic on~$\phi^{-1}(V)$ in~$M$.  
\endproclaim

It is easy to derive the partial differential equations that characterize 
harmonic morphisms in terms of local coframings.  Suppose that $U$ and $V$ 
are endowed with orthonormal coframings as above.  Suppose, further, that $v$ 
is any harmonic function on~$V$ and keep the notation as in the previous 
paragraph. Set $u = \phi^*(v) = \overline v$.  Pulling back the relation 
$dv = v_i\,\eta_i$ then yields
$$
du = f_{ia}\,\overline{v_i}\,\omega_a
\tag5
$$
so~$u_a = f_{ia}\,\overline{v_i}$.  Differentiating this relation and
comparing the result with the pullback of the relation~
$dv_i=v_k\,\eta_{ki}+v_{ij}\,\eta_j$ then yields
$$
u_{ab} = f_{ia}f_{jb}\,\overline{v_{ij}} + f_{iab}\,\overline{v_i}.
$$
In particular, 
$$
\Delta u = -f_{ia}f_{ja}\,\overline{v_{ij}} - f_{iaa}\,\overline{v_i}.
\tag6
$$
Now, because of the above-mentioned flexibility in choosing
the $2$-jet of a local harmonic function, it follows that a necessary
and sufficient condition for $\phi:U\to V$ to be a harmonic morphism 
is that there exist a function~$R$ on~$U$ so that the right hand side 
of~(6) is of the form $-R\,\overline{v_{ii}}$.  In other words,
$$
 f_{iaa} = 0
\qquad\text{and}\qquad
 f_{ia}f_{ja} = R\,\delta_{ij}
\quad\text{for some function~$R\ge0$ on~$U$.}
\tag7
$$
The first condition in~(7) is just that~$\phi$ be a harmonic map.
The second condition is known in the literature as {\it horizontal
weak conformality\/}, since it says that, for all~$p\in M$, the 
differential~$\phi'(p):T_pM\to T_{\phi(p)}N$ can be geometrically
described as orthogonal projection from the kernel of~$\phi'(p)$ to
its orthogonal complement followed by an isometry plus a conformal
scaling.   This characterization of harmonic morphisms can be found
in~\cite{EL1,(4.12)}, though it is due, independently, to Fuglede~\cite{Fu} 
and Ishihara~\cite{Is}.

Note that the conformal factor~$R$ defined in~(7) satisfies~$n\,R 
= ||\phi'||^2$ and hence vanishes only at the points where~$\phi'$ vanishes.  
According to~\cite{EL2, (2.32)}, if $R$ vanishes on an open set, 
then it vanishes identically on~$M$ (since $M$ is connected).  (This sort
of `unique continuation' principle holds more generally for all harmonic
maps, see~\cite{EL1, (3.16--18)}.)
Set aside this trivial case, in which~$\phi$ is constant, and
assume from now on that $\phi$ is non-constant. Thus, the set where~$R>0$ 
is a dense open subset~$M^*$ in~$M$.

By~(7), the harmonic morphism~$\phi:M^*\to N$ is a submersion, so
the dimension of~$M$ must satisfy~$m=n+p\ge n$.  The number~$p$
will be called the {\it fiber rank\/} or {\it corank\/} of~$\phi$. 

\subsubhead 1.1.3. 
Horizontal conformality
\endsubsubhead
Now, the conditions~(7) form an overdetermined system of~{\smc pde} for
the map~$\phi$.  The first condition, harmonicity, is second order and 
the second condition, horizontal (weak) conformality, is first order. 
The first step in studying this system is to examine the consequences
of the first order conditions.

Expand the summation convention to let lower case Greek indices run over 
the range~$n<\alpha,\beta,\gamma\le n{+}p$. (When $p=0$, this range is empty, 
so the formulae below have to be modified slightly in various obvious ways,
a task left to the reader.)

Assume that the open set~$V\subset N$ has a trivial
tangent bundle and that the rank~$p$ subbundle~$\ker\phi'\subset TU$ is 
trivial on the open set~$U\subset \phi^{-1}(V)\cap M^*$.
Let $r>0$ be the function on~$U$ satisfying~$n\,r^2=||\phi'||^2$.  
Then the 1-forms $\omega_i=r^{-1}\,\overline{\eta_i}=r^{-1}\phi^*(\eta_i)$
are $g$@-orthonormal on~$U$ and so, because of the triviality 
of the bundle~$\ker\phi'$, they can be completed to a $g$-orthonormal
coframing by choosing $p$ additional $1$@-forms~
$\omega_{n+1},\dots,\omega_{n+p}$
so that
$$
g_{|U} = \omega_1^2+\dots+\omega_{n+p}^2 = \omega_i^2 + \omega_\alpha^2\,.
$$
I will say that such a pair of coframings~$\omega$ on~$U$ and $\eta$ on~$V$
is {\it $\phi$@-adapted}. For a $\phi$@-adapted pair of coframings, (3)
simplifies to
$$
 f_{ia} = r\,\delta_{ia}.
\tag8
$$
  
Set~$r^{-1}\,dr = r_a\,\omega_a = r_i\,\omega_i+r_\alpha\,\omega_\alpha$.
Differentiating~$\overline{\eta_i} = r\,\omega_i$ yields
$$
\aligned
d\bigl(\overline{\eta_i}\bigr) & = dr\w\omega_i + r\,d\omega_i\\
r^{-1}\bigl(-\overline{\eta_{ij}}\w\overline{\eta_{j}}\bigr)
&=r^{-1}\,dr\w\omega_i-\omega_{ij}\w\omega_j-\omega_{i\alpha}\w\omega_\alpha\\
-\overline{\eta_{ij}}\w\omega_j
&=\bigl(r_i\,\omega_i{+}r_\alpha\,\omega_\alpha\bigr)\w\omega_i
-\omega_{ij}\w\omega_j-\omega_{i\alpha}\w\omega_\alpha\\
0&= \left[\bigl(r_k\,\omega_k{+}r_\alpha\,\omega_\alpha\bigr)\,\delta_{ij}
+\overline{\eta_{ij}}-\omega_{ij}\right]\w\omega_j
-\omega_{i\alpha}\w\omega_\alpha\,.\\
\endaligned
\tag9
$$
Since $\omega_{i\alpha}\w\omega_\alpha\equiv0\mod\omega_1,\dots,\omega_n$,
Cartan's Lemma implies that there exist functions~$A_{ij\alpha}$ and~
$H_{i\alpha\beta}=H_{i\beta\alpha}$ on~$U$ so that
$$
\omega_{i\alpha}=A_{ij\alpha}\,\omega_j + H_{i\alpha\beta}\,\omega_\beta\,.
\tag10
$$

The functions~$A_{ij\alpha}$ and $H_{i\alpha\beta}$ are 
the local components of globally defined tensor fields:  The expression
$$
{\bold H} = H_{i\alpha\beta}\,\omega_\alpha\circ\omega_\beta\otimes\ebold_i
$$
is independent of choice of frame field and has the defining property 
that it restricts to each fiber of~$\phi$ to be the second fundamental
 form of that fiber in~$M^*$.  The $ij$@-skewsymmetrization of the expression
$$
{\bold A} = A_{ij\alpha}\,\omega_i\otimes\omega_j\otimes\ebold_\alpha
$$
is the integrability tensor of the $n$-plane field~
$H_\phi=\bigl(\ker\phi'\bigr)^\perp$.  (The symmetric part will be
examined below.)

Substituting~$(10)$ into~$(9)$ yields the relation
$$
0 = 
\left[\bigl(r_k\,\omega_k{+}r_\alpha\,\omega_\alpha\bigr)\,\delta_{ij}
+\overline{\eta_{ij}}
-\omega_{ij}+A_{ij\alpha}\,\omega_\alpha\right]\w\omega_j
$$
and applying Cartan's Lemma to this yields that there exist functions
$s_{ijk}=s_{ikj}$ so that
$$
\bigl(r_k\,\omega_k{+}r_\alpha\,\omega_\alpha\bigr)\,\delta_{ij}
+\overline{\eta_{ij}}-\omega_{ij}+A_{ij\alpha}\,\omega_\alpha
= s_{ijk}\,\omega_k\,.
\tag11
$$
Symmetrizing the $\omega_\alpha$@-components on both sides of this
expression yields
$$
2\,r_\alpha + A_{ij\alpha} + A_{ji\alpha} = 0,
$$
implying that there are functions~$a_{ij\alpha}=-a_{ji\alpha}$ so that
$A_{ij\alpha} = -r_\alpha + a_{ij\alpha}$.  This gives a sharpened
version of~$(10)$, which now takes the form
$$
\omega_{i\alpha}
=-r_\alpha\,\omega_i+a_{ij\alpha}\,\omega_j+H_{i\alpha\beta}\,\omega_\beta\,,
\qquad\text{with $a_{ij\alpha}=-a_{ji\alpha}$.}
\tag$10'$
$$
Note that the skewsymmetric tensor~${\bold A}' = 
a_{ij\alpha}\,\omega_i\w\omega_j\otimes\ebold_\alpha$ is the 
integrability tensor for the $n$-plane field~$H_\phi$.

Substituting $(10')$ back into~$(11)$ allows one to solve for 
the ~$s_{ijk}$, finally yielding the relation
$$
\overline{\eta_{ij}} 
=\omega_{ij}+a_{ij\alpha}\,\omega_\alpha + r_j\,\omega_i-r_i\,\omega_j\,.
\tag12       
$$

\subsubhead 1.1.4. 
Harmonicity
\endsubsubhead
With this information derived from horizontal conformality, 
the quantities~$f_{iab}$ can now be computed via the defining formula~$(4)$:
$$
\align
 f_{iab}\,\omega_b 
&= df_{ia} -f_{ib}\,\omega_{ba} + f_{ja}\overline{\eta_{ij}}
 = d\bigl(r\delta_{ia}\bigr) - r\delta_{ib}\omega_{ba}
    +r\delta_{ja}\overline{\eta_{ij}}\\
&=r\left(\delta_{ia}\,\bigl(r_b\,\omega_b\bigr)
        -\omega_{ia}+\delta_{ja}\overline{\eta_{ij}}\right)\\
\endalign
$$
In particular, taking~$a=\alpha>n$ and using $(10')$,
$$
 f_{i\alpha b}\,\omega_b = -r\omega_{i\alpha}
=rr_\alpha\,\omega_i-ra_{ij\alpha}\,\omega_j-rH_{i\alpha\beta}\,\omega_\beta\,,
$$
so that~$f_{i\alpha\beta} = -rH_{i\alpha\beta}$.  Next, 
taking~$a=j\le n$, and using~$(12)$,
$$
\align
 f_{ijb}\,\omega_b 
&= r\left(\delta_{ij}\bigl(r_b\,\omega_b\bigr)
     -\omega_{ij}+\overline{\eta_{ij}}\right)\\
&= r\left(\delta_{ij}\bigl(r_b\,\omega_b\bigr)
     a_{ij\alpha}\,\omega_\alpha + r_j\,\omega_i-r_i\,\omega_j\right)\,.\\
\endalign
$$
so that, in particular, ~$f_{ijk} 
= r\bigl(\delta_{ij}r_k+\delta_{ik}r_j-\delta_{kj}r_i\bigr)$.

The components of the tension field~$\tau(\phi)$ can now expressed as
$$
 f_{iaa} = f_{ijj} + f_{i\alpha\alpha}
        = (2{-}n)\,rr_i - rH_{i\alpha\alpha}
        = -r\bigl((n{-}2)r_i+H_{i\alpha\alpha}\bigr).
\tag13
$$
Thus, a horizontally conformal submersion~$\phi$ is a harmonic morphism
if and only if it satisfies $(n{-}2)r_i+H_{i\alpha\alpha} = 0$.

\subsubhead 1.1.5. 
Elementary consequences
\endsubsubhead
Note that (13) yields the well-known result
that, when~$n=2$, a horizontally conformal map is a harmonic
morphism if and only if the fibers are minimal submanifolds, i.e., the
trace of the second fundamental form of each fiber vanishes~\cite{EL2,(2.34)}. 
Consequently, when~$n=2$, only the conformal structure on~$N$ is involved 
in determining whether or not the map~$\phi$ is a harmonic morphism, not
the full metric on~$N$.

Another consequence involves the case $p=0$, for, in this case, (13) has no 
terms~$H_{i\alpha\alpha}$.  Thus, in order for~$\phi$ to be a non-constant 
harmonic morphism between two manifolds of 
equal dimension~$n\not=2$, the terms~$r_i$ must all vanish, so that
the conformal factor~$r$ is constant.  Thus, when~$n\not=2$,
any non-constant harmonic morphism~$\phi:M^n\to N^n$ must be a homothety,
i.e., satisfy~$\phi^*h = r^2\,g$ for some constant~$r>0$.  

On the other hand, when~$p=0$ and $n=2$, equation~$(13)$ shows that
$f_{iaa}=0$ is an identity, so that the only condition for a map between 
Riemannian surfaces to be a harmonic morphism is that the map be a 
branched conformal map.  This is also well-known~\cite{EL1}.

Finally note that, when $n=1$, the condition of horizontal conformality
is automatic, so that a harmonic morphism is the same thing as a 
harmonic map.  

\remark{Remark}
One geometric interpretation of~(13) is that, for a harmonic morphism~
$\phi:M^*\to N$, parallel translation conserves volume in the fibers of~$\phi$
when it is measured with respect to the $p$@-form~
$r^{2-n}\,\omega_{n+1}\w\dots\w\omega_{n+p}$. This observation will be useful
in~\S2.
\endremark

\subhead 1.2. Moving frame computations for conformal foliations \endsubhead
The goal of this subsection is to derive necessary and sufficient conditions
on the invariants of a foliation of a Riemannian manifold~$M$ in order that 
its leaves be locally the fibers of some submersive harmonic morphism. 

It is useful to introduce a bit of terminology about foliations.  
A foliation~$\cF$ of codimension~$n$ on a manifold~$M^{n+p}$ is said to be
{\it amenable\/} if the leaf space~$M/\cF$ can be given the structure
of a smooth $n$@-manifold in such a way that the leaf projection
$M\to M/\cF$ becomes a smooth submersion.  When a smooth structure on~$M/\cF$ 
does exist with this property, it is unique.  Every foliation is locally
amenable in the sense that, for any smooth foliation~$\cF$ on~$M$, 
each point of~$M$ has an open neighborhood~$U$ on which~$\cF$ is amenable.
Such an open set~$U$ will be said to be {\it $\cF$@-amenable}.

Now, suppose that~$\cF$ be a foliation of codimension~$n$ on a manifold~
$M^{n+p}$ endowed with a Riemannian metric~$g$.  Using the metric~$g$,
the tangent bundle of~$M$ can be split as a $g$@-orthogonal direct sum
$$
TM = T\cF\oplus N\cF
$$
where~$T\cF$ and $N\cF$ are the tangent and normal bundles of the foliation~
$\cF$, respectively.  Correspondingly, the metric~$g$ can be split into 
two terms
$$
g = g' + g''
$$
where the quadratic form~$g'$ has $T\cF$ as its
null space and~$g''$ has ~$N\cF$ as its null space.

\subsubhead 1.2.1. 
Conformal foliations
\endsubsubhead
Let~$U\subset M$ be an open set on which~$T\cF$~and~$N\cF$ are both
trivial.  Then there exists a $g$-orthonormal coframing~
$\omega = (\omega_1,\dots,\omega_{n+p})$ on~$U$ with dual frame field~
$\ebold = (\ebold_1,\dots,\ebold_{n+p})$  so that the equations
$\omega_1 = \omega_2 = \dots = \omega_n = 0$ define~$T\cF$ while the
equations~$\omega_{n+1} = \omega_{n+2} = \dots = \omega_{n+p} = 0$ define~
$N\cF$.  In particular, the formula
$$
g'=\omega_1^2+\dots+\omega_n^2
$$
holds on~$U$. The equations~$d\omega_a = -\omega_{ab}\w\omega_b$ imply
(keeping the previously established index ranges for the summation
convention)
$$
d\omega_i = -\omega_{ia}\w\omega_a 
\equiv -\omega_{i\alpha}\w\omega_\alpha \mod\omega_1,\dots,\omega_n\,.
$$
Since the ~$\omega_i$ are the annihilators of a foliation, 
$d\omega_i\equiv0\mod\omega_1,\dots,\omega_n$. 
Thus, $\omega_{i\alpha}\w\omega_\alpha\equiv0\mod\omega_1,\dots,\omega_n$, 
so by Cartan's Lemma there exist functions~$A_{ij\alpha}$ 
and~$H_{i\alpha\beta}=H_{i\beta\alpha}$ so that
$$
\omega_{i\alpha} 
= A_{ij\alpha}\,\omega_j + H_{i\alpha\beta}\,\omega_\beta\,.
\tag1
$$

Now, the discussion from~\S1.1 shows that a necessary condition that the 
leaves of~$\cF$ in~$U$ be the fibers of a submersive harmonic morphism 
with domain~$U$ is that there exist functions~$r_\alpha$ on~$U$ so that
$$
A_{ij\alpha} + A_{ji\alpha} = -2r_\alpha\,\delta_{ij}\,.
\tag2
$$

Geometrically, condition~$(2)$ 
can be interpreted as the condition that the Lie derivative
of~$g'$ with respect to any vector field~$X = x_\alpha\,\ebold_\alpha$ 
(tangent to the leaves of~$\cF$) should be a multiple of~$g'$.  
Indeed, calculation yields
$$
\aligned
\Lie_X(g') &= 2 (\Lie_X\omega_i)\circ\omega_i
     = 2 (X\lefthook d\omega_i)\circ\omega_i
     = -2 \bigl(X\lefthook(\omega_{ia}\w\omega_a)\bigr)\circ\omega_i\,\\
&=
-2\bigl(\omega_{ia}(X)\,\omega_a-\omega_a(X)\,
\omega_{ia}\bigr)\circ\omega_i\,\\
&=
-2\bigl(\omega_{i\alpha}(X)\,\omega_\alpha
-x_\alpha\,\omega_{i\alpha}\bigr)\circ\omega_i\,\\
&=
-2\bigl(H_{i\alpha\beta}x_\beta\,\omega_\alpha
-x_\alpha\bigl(A_{ij\alpha}\,\omega_j 
+ H_{i\alpha\beta}\,\omega_\beta\bigr)\bigr)\circ\omega_i\,\\
&=x_\alpha(A_{ij\alpha}+A_{ji\alpha})\,\omega_i\circ\omega_j\,,\\
\endaligned
\tag3
$$
and this final expression is a multiple of $g'$ for any choice of the
 functions~$x_\alpha$ if and only if (2) is satisfied for some functions~
$r_\alpha$ on~$U$.  

This local condition can be expressed globally on~$M$.  If (2) is 
satisfied for some functions~$r_\alpha$, then the 1-form~$\rho''
=r_\alpha\,\omega_\alpha$ satisfies
$$
\Lie_X(g') = -2\rho''(X)\,g'
\tag$2'$
$$
 for all vector fields~$X$ tangent to the leaves of~$\cF$.  Now, the 1-form
$\rho''$ is independent of the choice of~$\omega$.  Thus, the global version
of~(2) is that there should exist a 1@-form~$\rho''$ on~$M$ (necessarily 
unique) that vanishes on~$N\cF$ and so that~$(2')$ is satisfied for all 
vector fields~$X$ tangent to the leaves of~$\cF$.  

Geometrically, condition~$(2')$ is satisfied if and only if, up to a 
conformal factor, the quadratic form~$g'$ can be pushed down onto the local 
leaf space~$U/\cF$ for any $\cF$@-amenable~$U\subset M$.  A foliation~$\cF$
satisfying~$(2')$ with respect to a given metric~$g$ is sometimes called a 
{\it conformal foliation\/} and this constitutes a first necessary condition 
in order for the leaves of~$\cF$ to be the fibers of a harmonic morphism.

\subsubhead 1.2.2. 
Local sufficient conditions
\endsubsubhead
To go further, the analysis must be divided into the cases~$n=2$ and~$n\not=2$.

 First, consider the case~$n=2$.  In~\S1.1, it was shown that
the non-singular fibers of a harmonic morphism~$\phi:M^{2+p}\to N^2$ would
be minimal submanifolds of codimension~$2$, so it is necessary
that the leaves of a codimension~$2$ foliation~$\cF$ on~$M$ be minimal
in order for them to be the fibers of a (local) submersive harmonic morphism.
Of course, $\cF$ must also be a conformal foliation.

Conversely, as has long been known~\cite{W1} (for example), these two 
conditions are locally sufficient.

\proclaim{Proposition 1}  Suppose that $\cF$ is a codimension~$2$ 
conformal foliation of~$M^{2+p}$ with the property that its leaves are minimal.
Then, for every open~$U\subset M$ that is~$\cF$@-amenable, there is a
conformal structure on the 2@-dimensional manifold~$U/\cF$ so that the 
leaf projection~$\phi: U\to U/\cF$ is a harmonic morphism.
\endproclaim

\demo{Proof}
This follows immediately from the discussion so far.   If~$U\subset M$ 
is an $\cF$@-amenable open set, then $(2')$ implies that $U/\cF$ carries a 
unique conformal structure so that the leaf projection~$\phi$ is horizontally
conformal.  From the discussion in~\S1.1 after equation~(13), the only other 
condition for~$\phi$ to be a harmonic morphism is that its fibers 
be minimal.\qed
\enddemo

Second, consider the case~$n\not=2$ and restrict to~$U$ with a coframing
as in~\S1.2.1.  Define functions~$r_i$ on~$U$ by~$(n{-}2)r_i= 
-H_{i\alpha\alpha}$.  The vector field~${\bold R} = r_i\,\ebold_i$ is, 
up to a constant scale factor, the mean curvature normal vector field for 
the leaves of~$\cF$.   Let~$\rho' = r_i\,\omega_i$.  This 1@-form is
locally defined with respect to a coframing, but, since it is the $g$@-dual
of the mean curvature vector field for the leaves of~$\cF$, it is, in fact,
globally defined on~$M$.  According to~(13) of~\S1.1, for any harmonic 
morphism~$\phi:M\to N$ whose fibers are the leaves of~$\cF$, the $1$@-form
$$
\rho = r_a\,\omega_a = r_i\,\omega_i + r_\alpha\,\omega_\alpha
     = \rho'+\rho'' 
     = \frac{-1}{(n{-}2)}\,H_{i\alpha\alpha}\,\omega_i 
        + r_\alpha\,\omega_\alpha
\tag4
$$
must satisfy~$\rho = d\bigl(\log r\bigr)$ where $r^2\,g' = \phi^*h$
and where $h$ is the metric on the target manifold~$N$.  Thus, a necessary
condition for the leaves of~$\cF$ to be the fibers of a harmonic morphism 
is that
$$
d\rho = 0.
\tag5
$$
The interesting result is the converse (which appears to be new):

\proclaim{Proposition 2}
If a foliation~$\cF$ of codimension~$n\not=2$ on a Riemannian manifold~
$(M,g)$ is conformal and satisfies~$(5)$, then for any $1$@-connected, 
$\cF$@-amenable open set~$U\subset M$, there is a metric~$h$ on~$U/\cF$, 
unique up to a constant scale factor, for which the leaf projection~
$\phi:U\to U/\cF$ is a harmonic morphism~$\phi:\bigl(U,g\bigr)
\to\bigl(U/\cF,h\bigr)$
\endproclaim

\demo{Proof}
Let~$U\subset M$ be 1@-connected and~$\cF$@-amenable and 
let~$\phi:U\to U/\cF$ be the leaf projection.  Since~$\cF$ is a 
conformal foliation, $\rho''$ can be defined so that $(2')$ holds while
$\rho''$ vanishes on~$N\cF$.  Moreover, let $\rho'$ be defined as
above to be the (suitably scaled) dual to the mean curvature vector field
of~$\cF$. Finally, use~(4) to define the 1@-form~$\rho$.  The assumption
that~$(5)$ holds is just that~$d\rho=0$.  

By the simple connectivity of~$U$, there exists a smooth positive function~$r$ 
on~$U$ so that~$\rho=d\bigl(\log r\bigr)$.  Since~$U$ is 
connected, $r$ is unique up to a multiplicative constant.

Let~$g'$ be the leaf normal part of the metric~$g$, as defined above. 
The formula~$(2')$ can be written as 
$$
\Lie_X(g') = -2\rho(X)\,g' = -2\,\Lie_X(\log r)\,g'
$$
 for all vector fields~$X$ tangent to the leaves of~$\cF$. It follows 
that~$\Lie_X (r^2\,g') = 0$ for all such vector fields~$X$. In particular, 
there must exist a metric~$h$ on~$U/\cF$ so that $r^2\,g' = \phi^*(h)$. 

The leaf projection~$\phi:\bigl(U,g\bigr)\to\bigl(U/\cF,h\bigr)$ is 
horizontally conformal by construction and equation~(13) of~\S1.1 now shows 
that $\phi$ is harmonic.  Thus, $\phi$ is a harmonic morphism.

The uniqueness of~$h$ with this property is immediate.
\qed
\enddemo  

By Propositions~1~and~2, the foliation~$\cF$ alone carries 
enough information to construct any harmonic morphism whose fibers are the
leaves of~$\cF$, locally and essentially uniquely, up to a constant 
scale factor on the range.  This observation will be important in the 
remainder of this article.  In particular, the classification of harmonic 
morphisms with a given domain Riemannian manifold~$(M,g)$ can be reduced to 
the classification of conformal foliations on domains in~$M$ that
either have codimension~2 and minimal leaves or else have codimension~
$n\not=2$ and satisfy~(5).

\head 2. Harmonic morphisms of corank~$1$: local theory  \endhead

The first result of this section will be Theorem~1, which says, roughly, 
that the local harmonic morphisms with a specified Riemannian $n$@-manifold~
$(N,h)$ as range depend on one arbitrary function of $n{+}1$ variables.  
The second result will be Theorem~2, which says, roughly, that, for
a specified Riemannian $(n{+}1)$@-manifold~$(M,g)$ with~$n\ge3$, the harmonic
morphisms with domain~$M$ form a finite dimensional space of dimension
at most ${{n+3}\choose2}-1$.  In~\S3, this space will be computed
explicitly when $g$ is a metric of constant sectional curvature on~$M$.

Since the assumption $p=1$ will hold for nearly all of the remainder of 
this article, it seems advisable to introduce a notational simplification: 
Instead of letting the index range run from~$1$ 
to~$n{+}1$, let it run from~$0$ to~$n$, i.e., make `$0$' an indicial 
synonym for~`$n{+}1$'.

\subhead 2.1. A local normal form \endsubhead
All harmonic morphisms of corank one can be put into a simple local
normal form.

\proclaim{Theorem 1}
Let~$\pi:P^{n+1}\to N^n$ be a principal $G$@-bundle with
$G = \bbR$ or $S^1$, let~$\psi$ be a connection form on~$P$, let
$h$ be a Riemannian metric on~$N$, and let~$r$ be a positive function
defined on an open set~$M\subset P$.  Then, when $M$ is endowed
with the metric~$g = r^{-2}\,\pi^*(h) + r^{2n-4}\,\psi^2$, the projection
$\pi:M\to N$ is a harmonic morphism.  Moreover, every submersive harmonic
morphism of corank~$1$ is locally of this form.
\endproclaim

\demo{Proof}
That $\pi:M\to N$ as described will be a harmonic morphism for any
choice of the data~$(h,\psi,r)$ is verified by a local calculation
using the formulae of~\S1.1.  I leave this to the reader, though
the reason that the formula works will become clearer after the 
converse part of the Proposition is argued.

It remains to prove that the given local normal form holds for
all submersive harmonic morphisms of corank~$1$.  Thus, 
let~$\phi:(M^{n+1},g)\to (N^n,h)$ be a submersion of Riemannian manifolds 
which is also a harmonic morphism.  Choose a point~$x\in M$ and let~
$V\subset N$ be a contractible open neighborhood of~$\phi(x)=y$. 
Let~$U\subset\phi^{-1}(V)\subset M$ be an open neighborhood of~$x$ on 
which the line bundle~$\ker\phi'$ is trivial and for
which the fibers of~$\phi:U\to V$ are connected and simply connected.  
Define $r>0$ by the equation $n\,r^2=||\phi'||^2$.

Let $\eta_1,\dots,\eta_n$ be an $h$-orthonormal coframing on~$V$ and let
$\omega_0$ on $U$ be chosen so that 
$$
\omega_0,
\ \omega_1 = r^{-1}\,\phi^*(\eta_1),
              \ \dots,\ \omega_n = r^{-1}\,\phi^*(\eta_n)
\tag1
$$
is a $g$-orthonormal coframing of~$U$, as was discussed in~\S1.1.  

Keep the notation of~\S1.1.3.  Since $\phi$ is a harmonic morphism,
by~(13) of that section,
$$
H_{i\alpha\alpha}=H_{i00} = -(n{-}2)\,r_i\,.
$$ 
Now, by equation~$(10')$ in~\S1.1, there exist functions~$a_{ij}=-a_{ji}$ 
(which would have been written as~$a_{ij0}$ in~\S1.1) 
and~$r_0,r_1,\dots,r_n$ so that
$$
\aligned
dr &= r\bigl(r_0\,\omega_0+r_1\,\omega_1+\dots+r_n\,\omega_n\bigr)\,,\\
\omega_{i0} &= - r_0\,\omega_i +a_{ij}\,\omega_j - (n{-}2)r_i\,\omega_0\,.\\
\endaligned
\tag2
$$
The structure equations then yield
$$
\aligned
d\bigl(r^{2-n}\,\omega_0\bigr)
&= (2{-}n) r^{-1}dr\w r^{2-n}\,\omega_0 + r^{2-n} d\omega_0\,,\\
&=r^{2-n}\left((2{-}n)(r_a\,\omega_a)\w
      \omega_0-\omega_{0i}\w\omega_i\right)\,,\\
&=r^{2-n}\left((2{-}n)(r_i\,\omega_i)\w\omega_0
  +\bigl(-r_0\,\omega_i
          +a_{ij}\,\omega_j-(n{-}2)r_i\,\omega_0\bigr)\w\omega_i\right)\,,\\
&=-r^{2-n}a_{ij}\,\omega_i\w\omega_j\,.\\
\endaligned
\tag3
$$ 
Set $\overline\Omega = -r^{2-n}a_{ij}\,\omega_i\w\omega_j$ and note that this 
is a closed $2$@-form which is $\phi$@-semibasic.  Since the fibers of~
$\phi$ in~$U$ are connected, there exists a (closed) 
$2$@-form~$\Omega$ on~$V$ so that~$\overline\Omega = \phi^*\Omega$.  Since
$V$ is contractible, there exists a $1$@-form~$\psi_0$
on~$V$ so that~$d\psi_0=\Omega$.  On $P = V\times\bbR$ (with coordinate~$t$ 
on the $\bbR$-factor), the $1$@-form~$\psi = dt+\psi_0$ represents
a connection for~$P$ regarded as a principal $\bbR$@-bundle over~$V$,
with projection onto the first factor~$\pi:P\to V$.  

Locally, $\bigl(U,r^{2-n}\,\omega_0\bigr)$ and $\bigl(P,\psi\bigr)$ are 
principal $\bbR$@-bundles over~$V$ with connection and they have the same 
curvature, so they are locally gauge equivalent.  Since~$P$ is globally
an~$\bbR$@-bundle over~$V$ and the fibers of~$\phi:U\to V$ are contractible, 
there exists a unique diffeomorphic embedding~$\tau:U\to P$ satisfying~
$\phi=\pi\circ\tau$ and~$\tau^*(\psi)=r^{2-n}\,\omega_0$ as well as
the (`gauge fixing') initial condition~$\tau(x) = \bigl(\phi(x),0\bigr)
\in V\times\bbR = P$.  Define~$s>0$ on~$\tau(U)\subset P$ so 
that~$\tau^*s = r$. Then
$$
g_{|U} = \tau^*\bigl(s^{-2}\,\pi^*(h) + s^{2n-4}\,\psi^2\bigr).
$$
Thus, the theorem is proved. \qed
\enddemo

\example{Example 1} {\it Isometric Quotients.}  An important special case
of this construction is the case where~$r$ is constant on the
 fibers of~$P$.  In this case, let ~$\Xbold$ be the ~$\pi$@-vertical
vector field that satisfies~$\psi(\Xbold) = 1$, i.e., $\Xbold$ is
the infinitesimal generator of the~$G$@-action.  Then, since~$r$
is constant on the flow lines of~$\Xbold$, it follows that~$\Xbold$
is a Killing vector field for~$g$.

Conversely, suppose that~$\Xbold$ is a Killing vector field 
on~$\bigl(M^{n+1},g\bigr)$, i.e., that $\Lie_\Xbold g = 0$.  
Let~$U\subset M$ be any open set on which~$\Xbold$ is non-vanishing and 
on which the foliation~$\cF_\Xbold$ by integral curves of~$\Xbold$ is 
amenable, with submersive leaf projection~$\pi:U\to N^n$. 

Let~$\psi$ be the $1$@-form that is ~$g$@-dual to~$\Xbold$ and let~$r>0$ 
be defined so that~$\omega_0 = r^{n-2}\,\psi$ is a unit $1$@-form. 
Because~$\Xbold$ is a symmetry vector field for~$g$, it follows that
$\Lie_\Xbold\psi=\Lie_\Xbold r = 0$, so in particular,~$\Lie_\Xbold\omega_0
=0$.  Moreover, since~$\psi(\Xbold)=1$, it follows that~$\psi$ is
a connection form for~$U$ regarded as a local principal~$\bbR$@-bundle
over~$N$.
 
Set~$g'= g - {\omega_0}^2$, as consistent with previous usage.
Since~$\Lie_\Xbold g = 0$, it follows that $\Lie_\Xbold g'=0$.
In particular, there exists a metric~$h$ on~$N^n$ so that $\pi^*(h) =
r^2 g'$.  (Note the use of~$\Lie_\Xbold r = 0$.)  Since, by construction,
$$
g = r^{-2}\,\pi^*(h) + r^{2n-4}\,\psi^2,
$$
it follows from Theorem~1 that the projection~$\pi:\bigl(U,g\bigr)\to
\bigl(N,h\bigr)$ is a harmonic morphism.

In this way, every non-zero Killing vector field defines a local harmonic
morphism of corank~1.  Note, however, that not every harmonic morphism
of corank~1 arises in this way.
\endexample

\remark{Remark: Local generality}
With Theorem~1 in hand, it can now be explained what is mean by
the statement, ``Up to local diffeomorphism, the harmonic morphisms
of corank~$1$ depend on one arbitrary function of $(n{+}1)$ variables.''

Consider the data~$(h,\psi,r)$ described in Theorem~1.  The first two 
components are locally defined on an $n$-manifold and thus, 
up to diffeomorphism, are specified by a certain number of functions of
$n$ variables.  However, $r>0$ is essentially arbitrary on an open set
of dimension~$(n{+}1)$.  This arbitrary function cannot be `gauged
away' by reparametrizing~$M$ since $\phi$ induces on~$M$ the local structure
of a principal $\bbR$@-bundle over an $n$@-dimensional base, and the
diffeomorphisms that preserve this local bundle structure depend only 
on functions of $n$ variables.

This sort of `dependency' discussion can be made more precise in
several ways.  One way is through the language of exterior
differential systems~\cite{BCG}, but it can also be interpreted in
terms of a formula for the dimension of the space of $k$@-jets of
metrics on a neighborhood of~$\bold0\in\bbR^{n+1}$ that admit a
submersive harmonic morphism to some $n$@-manifold, as in the theory
of Spencer, Goldschmidt, and Malgrange.  However, no essential use
of this notion of dependency or generality will be made in the rest of 
this article, so making it precise does not appear to be important.
Rather, this notion of local generality will be used in a heuristic
way for motivational purposes.
\endremark

\remark{Remark: Morphisms with a given range metric}
By Theorem~1, any Riemannian $n$@-manifold~$N$ can be the 
range of a submersive harmonic morphism of corank~1 in infinitely many 
distinct ways.  In fact, by judicious choice of~$r$, one can even
arrange that the domain~$\bigl(M^{n+1},g\bigr)$ be a complete
Riemannian manifold, regardless of whether~$(N,h)$ is complete or not.
\endremark

\remark{Remark: Higher corank}
There is a local structure theorem similar to Theorem~1 for submersive
harmonic morphisms of higher corank as well, but it is not quite as 
satisfactory.  

Using the `volume preserving' characterization mentioned
at the end of \S1.1, one can show that if~$\phi:(M^{n+p},g)\to(N^n,h)$ is a
harmonic morphism and ~$m\in M$ is a point where~$\phi'(m)\not=0$, then
there exists an $m$@-neighborhood~$U\subset M$ together with 
coordinates~$x_1,\ldots,x_{n+p}$ on~$U$ so that the following properties
hold:
\roster
\item The fibers of~$\phi$ in $U$ are defined by the equations~$dx_i=0$.
\item There exist functions $h_{ij}=h_{ji}$ of the variables~$x_i$ 
      so that 
	  $$ \phi^*(h) = h_{ij}\,dx_i\circ dx_j\,.
	  $$
\item There exist functions~$R>0$,~$P_{\alpha i}$, and~
      $g_{\alpha\beta}=g_{\beta\alpha}$ 
      of all the variables~$x_a$ satisfying
	  $$ \sum_\alpha\frac{\partial P_{\alpha i}}{\partial x_{\alpha}} = 0 
	      \quad\text{for all $i$,\quad and }\quad
		 \det(g_{\alpha\beta}) = 1
	  $$ 
	  so that
	  $$  g = R^{-p}\bigl(h_{ij}\,dx_i\circ dx_j\bigr) 
	      + R^{n-2} g_{\alpha\beta}\,
	       (dx_\alpha+P_{\alpha i}\,dx_i)\circ(dx_\beta+P_{\beta j}\,dx_j).
      $$
\endroster
The functions~$x_i$ can be regarded as coordinates on a neighborhood
of~$\phi(m)\in N$. In these coordinates, $\phi$ is just the 
natural submersion~$\bbR^{n+p}\to\bbR^n$ given by projection on the first
$n$ coordinates. 

Conversely, if $h_{ij}$, $R>0$, $P_{\alpha i}$, and $g_{\alpha\beta}$ are
 functions on some open set~$U\subset\bbR^{n+p}$ that satisfy the conditions 
in (2) and (3) and, moreover, the condition that the symmetric matrices~
$(h_{ij})$ and~$(g_{\alpha\beta})$ are positive definite at every point of~$U$,
then the metrics~$g$ on~$U$ and $h$ on~$V=\phi(U)\subset\bbR^n$ defined by
the above formulae have the property that~$\phi:(U,g)\to(V,h)$ is a 
harmonic morphism.

Note that the functions~$R$ and $g_{\alpha\beta}=g_{\beta\alpha}$ are 
`arbitrary', being only subject to algebraic inequalities and the single
algebraic equation~$\det(g_{\alpha\beta})=1$, while the $np$
 functions~$P_{\alpha i}$ are subject to $n$ linear, constant coefficient,
 first order {\smc pde}, but are otherwise arbitrary.

Thus, this normal form gives a local recipe for the general harmonic morphism~
$\phi:(M^{n+p},g)\to(N^n,h)$ that depends on the choice of
$$
1+\left({{p+1}\choose2}-1\right) + np - n = {{p+1}\choose2} + n(p-1)
$$
arbitrary functions of~$n{+}p$ variables.  However, it can be shown that
the ambiguity in choosing the `normal coordinates'~$x_a$ depends on~$p{-}1$
arbitrary functions of~$n{+}p$ variables.  Thus, this na{\"\i}ve counting
suggests that, up to local diffeomorphism, the general harmonic morphism
$\phi:(M^{n+p},g)\to(N^n,h)$ depends on
$$
{{p+1}\choose2} + (n-1)(p-1)
$$
arbitrary functions of~$n{+}p$ variables.  This count can be made rigorous,
but it will not be necessary to do so in this article.  The main point of 
these observations is that they indicate that, since, up to diffeomorphism, 
the general metric in dimension~$n{+}p$ depends on~${n+p}\choose2$ functions 
of $n{+}p$ variables rather than~${p+1}\choose2$ functions 
of $n{+}p$ variables, the general metric in this dimension will not admit, 
even locally, a harmonic morphism of corank~$p$ except when~$n=1$.

This normal form also indicates that, while any Riemannian manifold~$(N,h)$ 
can be the range of a harmonic morphism of arbitrary corank (in fact, in
infinitely many ways), the problem of finding the harmonic morphisms with a 
given Riemannian domain metric will involve the study of overdetermined 
systems and their integrability conditions.
\endremark

\subhead 2.2. Given domain metric \endsubhead
While Theorem~1 does give a local normal form for harmonic morphisms of
corank~1, it does not give much effective help for solving the problem 
of determining, for a specific Riemannian manifold~$\bigl(M^{n+1},g\bigr)$, 
what are the possible harmonic morphisms of the form~$\phi:M\to N$, where 
the $n$@-manifold~$N$ and its metric~$h$ may or may not be specified in 
advance.  This is essentially a question of when a given metric~$g$ can be 
written locally in the normal form given by Theorem~1 and, if so, in how many 
ways.

\subsubhead 2.2.1 
The low dimensional cases 
\endsubsubhead
The case~$n=1$ is basically trivial, and will not be discussed any
 further.

The case~$n=2$ is more interesting, but has a very different character
 from the cases where~$n\ge3$.  This is because, using the conformal
 flatness of metrics in dimension~2, the local normal form 
can be simplified when~$n=2$ to the coordinate form
$$
g=R(x_1,x_2,x_3)\,\bigl({dx_1}^2+{dx_2}^2\bigr)+
\bigl(dx_3+a(x_1,x_2)\,dx_1+b(x_1,x_2)\,dx_2\bigr)^2\,.
$$
where~$R>0$ is an arbitrary function of three variables while $a$ and $b$
are arbitrary functions of two variables and the fibers of the harmonic
morphism~$\phi$ are described by the equations~$dx_1=dx_2=0$.  

Computation using this form shows that, at any point in~$M$, 
the sectional curvatures of the $2$-planes passing through 
that point and containing the tangent to the $\phi$@-fiber must all be equal.  

However, elementary algebra shows that, if the sectional curvature is 
not constant on the $2$@-planes through a point~$x\in M$, 
then there are at most two tangent directions~$v\in T_xM$ 
so that the sectional curvature is constant on all the $2$@-planes 
containing~$v$.  Thus, for a metric~$g$ on~$M^3$ of non-constant sectional
curvature, there are at most two distinct foliations by curves which 
could be the fibers of a harmonic morphism of corank~1.  

It then suffices to apply the criteria developed in~\S1.2 to these 
two foliations (or one foliation when
there are only two distinct eigenvalues of the sectional curvature at
every point) to determine whether or not the given metric~$g$ admits a
harmonic morphism of corank~$1$.  This has been carried out by Baird and Wood~
\cite{BW2} in order to explicitly determine the metrics for which both 
of these foliations are the fibers of (possibly local) harmonic morphisms.

 Finally, when~$g$ has constant sectional curvature, there are
many harmonic morphisms~$\phi:U\to N^2$ of corank~$1$ with domain an open 
subset~$U\subset M$.  The local description is very simple:  
The fibers of such a map~$\phi$ must be geodesics (since they are
minimal) and so~$N^2$ can be identified with a surface~$\Sigma_\phi$ 
in the space~$\Lambda\bigl(M\bigr)$ of geodesics in~$M$.  
If~$M$ is complete and simply connected, this latter 
space has dimension~$4$ and carries a canonical complex structure which 
makes it into a complex surface.  The surface~$\Sigma_\phi$ is then seen
to be a complex curve in this complex surface.  Conversely, if $C\subset
\Lambda\bigl(M\bigr)$ is any complex curve (aside from a few special
cases), then the 2-parameter family of geodesics that it represents will
cover an open set in~$M$.  By suitably restricting the curve~$C$, one
can arrive at an open set~$U\subset M$ that is foliated by these geodesics
and this foliation will be amenable and satisfy the conditions of~\S1.2, 
so that it gives rise to a harmonic morphism with domain~$U$ and of corank~$1$.
 For more details and a study of the global and singularity issues, 
see~\cite{BW1}.

\subsubhead 2.2.2 
The high dimensional cases 
\endsubsubhead
The rest of this article concerns only the case~$n\ge3$.  The main
result to conclude this section will be a general `finiteness'
theorem, showing that the set of harmonic morphisms of corank~1 with
a given Riemannian domain~$(M^{n+1},g)$ is a finite dimensional space.

Let~$(M,g)$ be a Riemannian $(n{+}1)$@-manifold.  By Proposition~2, for any 
connected open set~$U\subset M$ on which there exists a 
submersive harmonic morphism~$\phi:(U,g)\to(N^n,h)$, the $n$@-manifold~$N$
and, up to a constant scale factor, the metric~$h$ are determined by 
the foliation~$\cF_\phi$ defined by the fibers of~$\phi$.  For simplicity, I
am going to assume that it is possible smoothly to orient the fibers 
of~$\phi$, as can always be arranged by passing to a double cover of~$U$
if necessary.  Then there will exist a unique unit vector field~$\ubold_\phi$
on~$U$ tangent to the fibers of~$\phi$ and inducing the given orientation 
on each fiber.

Conversely, given a unit vector field~$\ubold$ on an open set~$U\subset M$,
this process can be reversed.  Provided that~$\ubold$ satisfies
certain equations (to be spelled out below) and provided that the foliation 
$\cF_\ubold$ of~$U$ induced by the integral curves of~$\ubold$ is amenable, 
then, by Proposition~2, there will be a metric~$h$, unique up to a constant
multiple, on~$U/\cF_\ubold$ so that the leaf projection
$\phi:U\to U/\cF_\ubold$ is a harmonic morphism.

The conditions on~$\ubold$ needed to construct this harmonic morphism are
as follows:  Let~$\omega_0$ be the 1@-form on~$U$ that is $g$@-dual 
to~$\ubold$.  Set~$g' = g - {\omega_0}^2$.  Then the first condition on~
$\ubold$ is that the foliation~$\cF_\ubold$ should be conformal, i.e., that
there should exist a function~$r_0$ on~$U$ so that
$$
\Lie_\ubold(g') = -2r_0\,g'\,.
\tag4
$$
Note that~(4) is an overdetermined set of first order {\smc pde} for~$\ubold$.
The 1@-form~$\rho'' = r_0\,\omega_0$ is the same as the one constructed
in~\S1.2.1. Next, define a 1-form~$\rho'$ on~$U$ by the equation 
$$
\Lie_\ubold \omega_0 = (n-2)\,\rho'\,.
$$
Since $n\not=2$, this can always be done.  This~$\rho'$ is the same as
the one constructed in~\S1.2.2.  

The final condition, corresponding to~(5) of~\S1.2.2, is
$$
d\bigl( r_0\,\omega_0 + \rho' \bigr) = 0.
\tag5
$$
Note that~(5) is an overdetermined set of second order {\smc pde} for~$\ubold$.  

Conversely, by Proposition~2, conditions (4) and (5) are sufficient to 
imply that, locally, the integral curves of~$\ubold$ are the fibers of
a harmonic morphism.

To state the next theorem precisely requires the use of the language of
jets.  Let~$SM\to M$ denote the unit sphere bundle of~$M$.  This is a
bundle of fiber rank~$n$ over~$M$.  

Let~$J^1(M,SM)$ denote the bundle of 1@-jets of sections of~$SM$, 
and let~$\Sigma_1\subset J^1(M,SM)$ denote the subset consisting of those 
1@-jets of sections that satisfy~(4).  (This makes sense since (4) is a set 
of first order equations for a section~$\ubold$ of~$SM$.)  Then~$\Sigma_1$ 
is a smooth manifold of codimension~${{n+1}\choose2}-1$ in~$J^1(M,SM)$ 
and the basepoint projection~$\Sigma_1\to M$ makes~$\Sigma_1$ into a 
smooth bundle over~$M$ with fibers of dimension~${{n+2}\choose2}$.  

 Finally, let~$J^2(M,SM)$ denote the bundle of 2@-jets of sections of~$SM$, 
and let~$\Sigma_2\subset J^2(M,SM)$ denote the subset consisting of 
those 2@-jets of sections that satisfy~(4)~and~(5).  As will come out
in the proof of Theorem~2 below, $\Sigma_2$ is a smooth manifold
and the natural projection~$\Sigma_2\to\Sigma_1$ is a smooth submersion
with fibers diffeomorphic to~$\bbR^{n+1}$.  Thus, $\Sigma_2$ is a smooth 
bundle over~$M$ with fibers of dimension
$$
(n{+}1) + {{n+2}\choose2} = {{n+3}\choose2}-1\,.
$$

\proclaim{Theorem 2}
The bundle~$\Sigma_2\to M$ has a horizontal plane field~$H$ of 
dimension~$n{+}1$ with the property that, for every unit vector field~
$\ubold$ satisfying~$(4)$ and $(5)$, the section~$j^2\ubold$ of~$\Sigma_2$ 
is tangent to~$H$.  

In particular, if, on a connected open set~$U\subset M$, there are
two unit vector fields~$\ubold_1$~and~$\ubold_2$, each of which
satisfies~$(4)$~and~$(5)$, that agree to second order at some point of~$U$,
then $\ubold_1=\ubold_2$ throughout~$U$.
\endproclaim

\demo{Proof} Let~$F\to M$ be the $g$@-orthonormal frame bundle of~$M$.  The
elements of~$F$ are of the form~$(m;u)=(m;u_0,u_1,\ldots,u_n)$ where ~$m$ is a 
point of~$M$ and~$u=(u_0,u_1,\ldots,u_n)$ is a $g$@-orthonormal basis 
of~$T_mM$.  Then~$\pi:F\to M$ is a principal right~$\Or(n{+}1)$@-bundle 
over~$M$, where~$\pi$ is the basepoint projection and the right action 
by~$A=(A_{ab})$ in~$\Or(n{+}1)$ is the usual one:
$$
(m;\,u\,)\cdot A = (m;\,u_b\,A_{ba}\,).
$$
There are canonical 1-forms~$\omega_a$ on~$F$ that satisfy
$$
\omega_a(v) = \pi'(v)\cdot u_a
$$
 for every vector~$v$ in the tangent space to~$F$ at~$(m;u)$.  There
are also the Levi-Civita connection forms~$\omega_{ab}=-\omega_{ba}$ that
satisfy the first structure equations
$$
d\omega_a = -\omega_{ab}\w\omega_b
$$
and the second structure equations
$$
\Omega_{ab} = d\omega_{ab} + \omega_{ac}\w\omega_{cb} 
= {\ts\frac12}R_{abcd}\,\omega_c\w\omega_d
$$
where the functions~$R_{abcd}=-R_{bacd}=-R_{abdc}$ are well-defined on~$F$
and represent components of the Riemann curvature tensor in the sense 
that
$$
\pi^*\bigl(\text{Riem}(g)\bigr) = 
R_{abcd}\,\omega_a\otimes\omega_b\otimes\omega_c\otimes\omega_d\,.
$$
This use of the symbols~$\omega_a$~and~$\omega_{ab}$ should not be confused
with the earlier usages.  Previously, the forms~$\omega_a$~and~$\omega_{ab}$
were defined on an open set~$U\subset M$ and relative to a coframing, i.e., 
a section of~$F$ over~$U$.  

The map~$\sigma:F\to SM$ defined by~$\sigma(m;u)=(m,u_0)$ is a submersion
and makes~$F$ into a principal right $\Or(n)$@-bundle over~$SM$.

Let~$\ubold$ be a unit vector field on~$U\subset M$ and let $F_\ubold
=\sigma^{-1}\bigl(\ubold(U)\bigr)\subset F$.  Then $\pi:F_\ubold\to U$ is a 
principal right~$\Or(n)$@-bundle.  On~$F_\ubold$, the 1@-form~$\omega_0$
equals~$\pi^*\bigl(\ubold^\flat\bigr)$, where $\ubold^\flat$ is the 1@-form
that is $g$@-dual to~$\ubold$.  Therefore, on~$F_\ubold$, the closed 
2@-form~$d\omega_0$ must be $\pi$@-basic.  Thus, since
$$
d\omega_0 = -\omega_{0i}\w\omega_i = \omega_{i0}\w\omega_i \,,
$$ 
there must be functions~$A_{ij}$ and~$r_i$ on~$F_\ubold$ so that 
$$
\omega_{i0} = -(n{-}2)r_i\,\omega_0 + A_{ij}\,\omega_j\,.
$$
(The coefficient~$-(n{-}2)$ is introduced to align this formula
with previously derived ones.)  

Now~$g=\bigl(\ubold^\flat\bigr)^2+g'$ and calculation as in~(3) of~\S1.2.1
shows that~$\Lie_\ubold(g') = -2r_0\,g'$ if and only if the functions~
$A_{ij}$ can be written in the form
$$
A_{ij} = -r_0\,\delta_{ij} + a_{ij}
$$
where~$a_{ij} = -a_{ji}$.  Thus, for any unit vector field~$\ubold$, 
condition~(4) is equivalent to the condition that equations of the form
$$
\omega_{i0} = -(n{-}2)r_i\,\omega_0 + (-r_0\,\delta_{ij} + a_{ij})\,\omega_j
$$
hold on~$F_\ubold$ for some functions~$r_a$ and~$a_{ij}=-a_{ji}$.

Because the formula~$\pi^*\bigl(\nabla\ubold\bigr) = \omega_{i0}\otimes u_i$
holds, the functions~$r_a$ and~$a_{ij}$ represent the
`undetermined' derivatives of any solution~$\ubold$ of~(4).
This suggests that the submanifold~$\Sigma_1\subset J^1(M,SM)$ and 
its canonical contact ideal can be described by the following construction:  

Regard $\bbR^n$ as the set of 
columns of real numbers of height~$n$ and~$\euso(n)$ as
the space of skew-symmetric matrices of size~$n$@-by@-$n$.  Let
$$
X = F\times\bbR\times\bbR^n\times\euso(n)
$$ 
and let~$\rbold_0:X\to\bbR$, $\rbold = (\rbold_i):X\to\bbR^n$, 
and~$\abold = (\abold_{ij}):F\to\euso(n)$ denote the 
projections onto the second third and fourth factors of this product, 
respectively.

Define 1@-forms~$\thetabold_i$ on~$X$ by the formulae
$$
\thetabold_i
  =\omega_{i0}+(n{-}2)\rbold_i\,\omega_0
    +(\rbold_0\,\delta_{ij}-\abold_{ij})\,\omega_j\,.
$$
Let~$\cI_0$ be the differential ideal on~$X$ generated by the 1@-forms~
$\thetabold_i$.

 For vector fields~$\ubold$ that satisfy condition~(4), the bundle~$F_\ubold
\subset F$ can be lifted up to~$X$ by the mapping
$$
(m;u) \mapsto \bigl((m;u),r_0,(r_i),(a_{ij})\bigr)
$$
and, by construction, this lifting~$F_\ubold\hookrightarrow X$ is an 
integral manifold of~$\cI_0$ on which the forms~$\omega_a$ and~$\omega_{ij}$ 
are linearly independent.

Conversely, any connected integral manifold in~$X$ of the~$\thetabold_i$ on 
which the forms~$\omega_a$ and~$\omega_{ij}$ are linearly independent is the 
lift of an open subset of~$F_\ubold$ for some unit vector field~$\ubold$ 
satisfying~(4) that is defined (though possibly multi-valued) on some open 
subset of~$M$.  This can be verified directly, but, in any case, is a
standard argument in the theory of exterior differential systems. For
 further details, see~\cite{BCG}.

Let $\Or(n)$ act freely on the right on~$X$ by the formula
$$
x\cdot A = \bigl((m;u),\rbold_0,\rbold,\abold\bigr)\cdot A
= \bigl((m;u)\cdot A,\,\rbold_0,\,A^{-1}\rbold,A^{-1}\abold A\bigr).
$$ 
 for all~$A\in\Or(n)$. The $\Or(n)$@-orbits on~$X$ are the Cauchy 
characteristics of the Pfaffian system~$\cI_0$~\cite{BCG, Chapter 1}.  
Writing~$\thetabold=(\thetabold_i)$, the formula~
$R_A^*\thetabold=A^{-1}\,\thetabold$ holds for~$A\in\Or(n)$ 
and the form~$\thetabold$ is semibasic for the quotient projection~
$X\to X/\Or(n)$.  Thus, there is a well-defined Pfaffian system~
$\bar\cI_0$ on~$X/\Or(n)$ of rank~$n$ that pulls back to~$X$ to be~$\cI_0$.
It is important to understand the geometric meaning of the quotient~$X/\Or(n)$.

Recall that, for any smooth bundle~$B\to M$, the space~$J^1(M,B)$ of 1-jets
of sections of~$B$ can be identified with the space of pairs of the form~
$(b,P)$ where $b$ is an element of~$B$ and $P\subset T_bB$ is an
$(n{+}1)$@-plane transverse to the fiber of~$B\to M$ containing~$b$.  

Accordingly, a map~$X\to J^1(M,SM)$ can be defined:  For every~$x\in X$, 
let $\Pbold_x\subset T_xX$ be 
the codimension~$n$ subspace annihilated by the components of~$\thetabold$.
The subspace~$\Pbold_x$ contains the tangents to the $\Or(n)$@-orbits
and so pushes down to define a codimension~$n$ plane field~$\bar\Pbold$
on~$X/\Or(n)$ that is the annihilator of the 1@-forms in~$\bar\cI_0$.

Now, $\Pbold_x$ contains the tangents to the fibers of the submersion~
$\sigma_0:X\to SM$ defined by~$\sigma_0\bigl((m;u),r_a,a_{ij}\bigr)=(m,u_0)$, 
so, for every~$x\in X$, there is a well-defined subspace~$\sigma_1(x)
\subset T_{\sigma_0(x)}SM$ that is of codimension~$n$ and is transverse 
to the fibers of~$SM\to M$.  By the above identification,
this defines a smooth map~$\sigma_1:X\to J^1(M,SM)$.  Tracing through
the definitions, one finds that the image of~$\sigma_1$ is~$\Sigma_1$ while
the fibers of~$\sigma_1$ are the $\Or(n)$@-orbits in~$X$.  Thus,
$$
X/\Or(n) = \Sigma_1 \subset J^1(M,SM)
$$
and, under this identification, $\bar\cI_0$ is just the pullback
of the canonical contact system on~$J^1(M,SM)$ to~$\Sigma_1$.  This
identification explains the significance of the system~$\cI_0$.

Now, condition~(5) on the vector field~$\ubold$ corresponds to the condition
$$
d\bigl(r_0\,\omega_0 + r_i\,\omega_i\bigr) = 0
$$
on~$F_\ubold$.  This suggests defining a 2@-form~$\Thetabold_0$ on~$X$ 
by the formula
$$
\Thetabold_0 = d\bigl(\rbold_0\,\omega_0 + \rbold_i\,\omega_i\bigr).
$$
Note that~$\rbold_0\,\omega_0 + \rbold_i\,\omega_i$ is semibasic for the
projection~$X\to X/\Or(n)$ and invariant under the action of~$\Or(n)$. 
Thus, this 1@-form is the $\sigma_1$@-pullback of a 1@-form on~$\Sigma_1 
= X/\Or(n)$. In particular, its exterior derivative~$\Thetabold_0$ 
is the $\sigma_1$@-pullback of a 2@-form on~$\Sigma_1$.

Let~$\cI$ be the differential ideal on~$X$ generated by the (closed) 2@-form~
$\Thetabold_0$ and the 1@-forms~$\thetabold_i$.  As an algebraic
ideal, $\cI$ is generated by the 1-forms~$\thetabold_i$ and the 2@-forms~
$\Thetabold_0$ and $d\thetabold_i$.  Like~$\cI_0$, the Cauchy characteristics 
of~$\cI$ are the orbits of the $\Or(n)$ action on~$X$.  Thus,~$\cI$ is the 
$\sigma_1$@-pullback of a differential ideal~$\bar\cI$ on~$\Sigma_1$.  

By construction,  
the integral manifolds of~$\bar\cI$ that satisfy the independence condition~
$\Omega_1 = \omega_0\w\omega_1\w\dots\w\omega_n\not=0$ are locally the 
$1$@-jet graphs in~$\Sigma_1\subset J^1(M,SM)$ of unit vector fields~$\ubold$ 
on open domains in~$M$ that satisfy~(4)~and~(5).  Thus, it is the 
exterior differential system~$\bigl(\bar\cI,\Omega_1\bigr)$ that must be
studied in order to determine the space of solutions to~(4) and (5). 

 For the rest of the proof, I will be applying the standard exterior 
differential systems techniques to~$\bigl(\bar\cI,\Omega_1\bigr)$.  The
reader who wants to know more about these techniques can consult~\cite{BCG,
Chapter~4}.  In particular, by these methods, the space~$\Sigma_2$, i.e., 
the space of 2@-jets of solutions to the equations~(4) and~(5), is seen to be 
the space of integral elements of~$\bigl(\bar\cI,\Omega_1\bigr)$ on~$\Sigma_1$.
(Recall that, given an exterior differential system~$\cJ$ on a manifold~$M$
with an independence condition defined by some~$p$@-form~$\Omega$ on~$M$,
the  {\it integral elements\/} of~$\bigl(\cJ,\Omega\bigr)$ are the 
$p$@-planes~$E\subset T_mM$ on which the elements of~$\cJ$ vanish but on
which~$\Omega$ does not.)

It is more convenient to calculate on~$X$ than on~$X/\Or(n)$.  
Because the Cauchy characteristics of~$\cI$ are swept out by the 
$\Or(n)$@-action on~$X$, the integral manifolds of~$\bigl(\bar\cI,
\Omega_1\bigr)$ on~$\Sigma_1$ are in one@-to@-one correspondence with the 
integral manifolds of the system~$\bigl(\cI,\Omega_2\bigr)$ on~$X$, where
$$
\Omega_2= 
\Omega_1\w \bigl(\omega_{12}\w\omega_{13}\w\dots\w\omega_{(n-1)n}\bigr).
$$
In particular, $\Sigma_2$ can be computed from the structure equations 
of~$\cI$.

A routine, albeit tedious, calculation yields structure equations
$$
\Thetabold_0 \equiv \rhobold_0\w\omega_0 + \rhobold_i\w\omega_i
\mod \thetabold\,,
\tag6
$$ 
and 
$$
d\thetabold_i\equiv 
(n{-}2)\rhobold_i\w\omega_0
    +(\delta_{ij}\,\rhobold_0-\alphabold_{ij})\w\omega_j
\mod \thetabold\,,
\tag7
$$
where 
$$
\aligned
\rhobold_0 &= d\rbold_0\,,\\
\rhobold_i &= d\rbold_i - \rbold_j\,\omega_{ji}
               -\bigl((n{-}1)\,\rbold_0\rbold_i
                          +\abold_{ij}\rbold_j\bigr)\,\omega_0 
               +(\rbold_0\,\abold_{ij}-\rbold_{ij})\,\omega_j\,,\\
\alphabold_{ij} &= d\abold_{ij}
                   - \abold_{kj}\,\omega_{ki}-\abold_{ik}\,\omega_{kj}
                   - n\,\rbold_0\abold_{ij}\,\omega_0
                   - \bbold_{ijk}\,\omega_k\\
\endaligned
\tag8
$$
and where the expressions~$\rbold_{ij}=\rbold_{ji}$ and $\bbold_{ijk}$ 
are defined on~$X$ by the equations
$$
\bbold_{ijk}= (n{-}2)\,\bigl(\rbold_i\abold_{jk}+\rbold_j\abold_{ki}
                     +\rbold_k\abold_{ji}\bigr)
               +{\ts\frac12}\bigl(R_{0ijk}-R_{0jik}-R_{0kij}\bigr)
$$
and
$$
-(n{-}2)\,\rbold_{ij} 
= \delta_{ij}\,{\rbold_0}^2
      +\abold_{ik}\abold_{kj} +(n{-}2)^2\,\rbold_i\rbold_j + R_{i0j0}\,.       
$$

 From (6) and (7), any integral element~$E\subset T_xX$ of~
$\bigl(\cI,\Omega_2\bigr)$ is defined by equations
$$
\aligned
   \thetabold_i &= 0\,,\\
     \rhobold_0 &= {} - (n{-}2)\,s_0\,\omega_0 - s_i\,\omega_i\,,\\
     \rhobold_i &= {} - s_i\,\omega_0 - s_0\,\omega_i\,,\\
\alphabold_{ij} &= {} - s_i\,\omega_j + s_j\,\omega_i\\
\endaligned
\tag9
$$
 for some unique numbers~$s_a$ (that depend on~$E$).  Conversely, for any 
choice of numbers~$s_a$, equations~(9) define an integral element
of~$\bigl(\cI,\Omega_2\bigr)$ at every point of~$X$.

One can now compute the Cartan characters of~$\cI$ from structure 
equations~(6,7) and conclude that~$\cI$ is not involutive.  Consequently,
it will be necessary to examine the first prolongation of~
$\bigl(\cI,\Omega_2\bigr)$, i.e., the space~$Y$ of integral
elements of~$\bigl(\cI,\Omega_2\bigr)$ and its canonical differential
ideal~$\cI^{(1)}$. Now, by equations~(9), the space $Y$ is diffeomorphic to
$$
Y = X\times\bbR\times\bbR^n.
$$
Let~$\sbold_0:Y\to\bbR$ and $\sbold = (\sbold_i):Y\to\bbR^n$ be the
projections on the second and third factors respectively.  The first
prolongation ideal of~$\cI$ is then the ideal~$\cI^{(1)}$ on~$Y$ generated 
by the~$\thetabold_i$ and the $1$@-forms~$\etabold_a$ and $\etabold_{ij}=
-\etabold_{ji}$ defined by 
$$
\aligned
   \etabold_0&=\rhobold_0+(n{-}2)\,\sbold_0\,\omega_0+\sbold_i\,\omega_i\,,\\
   \etabold_i&=\rhobold_i+\sbold_i\,\omega_0+\sbold_0\,\omega_i\,,\\
\etabold_{ij}&=\alphabold_{ij}+\sbold_i\,\omega_j-\sbold_j\,\omega_i\,.\\
\endaligned
\tag10
$$

If~$N\subset X$ is any integral manifold of~$\bigl(\cI,\Omega_2\bigr)$,
there will be unique functions~$s_a$ on~$N$ so that equations~(9) hold 
(since the tangent spaces to~$N$ must be integral elements of~
$\bigl(\cI,\Omega_2\bigr)$). Thus, the lifting~$N\hookrightarrow Y$
defined by~$x\mapsto\bigl(x,s_0(x),s_i(x)\bigr)$ for~$x\in N$ lifts~$N$
to an integral manifold of~$\cI^{(1)}$.  Conversely, every integral manifold
of~$\bigl(\cI^{(1)},\Omega_2\bigr)$ on~$Y$ is the lift of a unique
integral manifold of~$\bigl(\cI,\Omega_2\bigr)$ on~$X$.  Because of this
one@-to@-one correspondence of integral manifolds, it suffices to determine 
the integral manifolds on~$Y$.

Now, the $\Or(n)$@-action on~$X$ lifts in the obvious way to an 
$\Or(n)$@-action on~$Y$ whose orbits are the Cauchy characteristics
of the system~$\cI^{(1)}$.  The quotient space~$Y/\Or(n)$ is naturally
identified with the set of integral elements of~$\bigl(\bar\cI,\Omega_1\bigr)$,
i.e., with~$\Sigma_2\subset J^2(M,SM)$.  By its very construction, the
second order contact system on~$J^2(M,SM)$ pulls back to~$\Sigma_2$ to
become an exterior differential system~$\bar\cI^{(1)}$ which, in turn, 
then pulls back to~$Y$ to become~$\cI^{(1)}$.  

The next step is to determine the structure equations of~$\cI^{(1)}$.
Now, equations~(6) and (7) can be written as
$$
\Thetabold_0 \equiv \etabold_0\w\omega_0 + \etabold_i\w\omega_i
\mod \thetabold\,,
\tag11
$$ 
and 
$$
d\thetabold_i\equiv 
(n{-}2)\etabold_i\w\omega_0
    +(\delta_{ij}\,\etabold_0-\etabold_{ij})\w\omega_j
\mod \thetabold\,.
\tag12
$$
Thus, $\Thetabold_0\equiv d\thetabold_i\equiv 0 \mod \thetabold,\etabold$.
To complete the structure equations of~$\cI^{(1)}$, formulae for
$d\etabold_a$ and $d\etabold_{ij}$ must now be computed.  

Since~$d\Thetabold_0=0$, differentiating (11) and reducing modulo
$\thetabold,\etabold$ yields
$$
d\etabold_0\w\omega_0 + d\etabold_i\w\omega_i
\equiv 0 \mod \thetabold,\etabold\,.
\tag13
$$
It follows that there exist 1-forms~$\sigmabold_a$ so
that
$$
d\etabold_0\equiv (n{-}2)\,\sigmabold_0\w\omega_0+\sigmabold_i\w\omega_i
 \mod \thetabold,\etabold\,,
\tag14
$$
and examination of the definition of~$\etabold_0$ reveals that
$$
\sigmabold_a\equiv d\sbold_a \mod \omega_b,\omega_{ij},\thetabold,\etabold\,.
\tag15
$$
(A more explicit formula for~$\sigmabold_a$ will not be needed in the proof.)
Substituting~(14) into~(13) and collecting terms yields
$$
\bigl(d\etabold_i-\sigmabold_i\w\omega_0\bigr)\w\omega_i
\equiv 0 \mod \thetabold,\etabold\,,
$$
which implies that there exist 1@-forms~$\tau_{ij}=\tau_{ji}$ so that
$$
d\etabold_i\equiv\sigmabold_i\w\omega_0 +\tau_{ij}\w\omega_j
\mod \thetabold,\etabold\,.
\tag16
$$
Now, differentiating~(12) and reducing modulo~$\thetabold$ and $\etabold$
yields
$$
0\equiv 
(n{-}2)d\etabold_i\w\omega_0
    +(\delta_{ij}\,d\etabold_0-d\etabold_{ij})\w\omega_j
\mod \thetabold,\etabold\,.
\tag17
$$
Substituting (14) and (16) into (17) and rearranging then yields
$$
0\equiv 
(n{-}2)\,\bigl(\tau_{ij}-\delta_{ij}\sigmabold_0\bigr)\w\omega_j\w\omega_0
 -\bigl(d\etabold_{ij}-\sigmabold_i\w\omega_j+\sigmabold_j\w\omega_i
             \bigr)\w\omega_j
\mod \thetabold,\etabold\,.
\tag18
$$
Reducing~(18) modulo~$\omega_0$ gives
$$
0\equiv \bigl(
            d\etabold_{ij}-\sigmabold_i\w\omega_j+\sigmabold_j\w\omega_i
             \bigr)\w\omega_j
\mod \thetabold,\etabold,\omega_0\,.
$$
Due to the skewsymmetry in~$i$ and~$j$ of the expression~
$d\etabold_{ij}-\sigmabold_i\w\omega_j+\sigmabold_j\w\omega_i$, it follows
that there are functions~$F_{ijkl}$ with the symmetries of a Riemann
curvature tensor so that
$$
d\etabold_{ij}\equiv \sigmabold_i\w\omega_j-\sigmabold_j\w\omega_i
                      +{\ts\frac12}F_{ijkl}\,\omega_k\w\omega_l
\mod \thetabold,\etabold,\omega_0\,.
$$
Thus, there exist 1-forms~$\psi_{ij}=-\psi_{ji}$ so that
$$
d\etabold_{ij}\equiv \sigmabold_i\w\omega_j-\sigmabold_j\w\omega_i
                      +(n{-}2)\,\psi_{ij}\w\omega_0
                      +{\ts\frac12}F_{ijkl}\,\omega_k\w\omega_l
\mod \thetabold,\etabold\,.
\tag19
$$
Substituting~(19) into (18) yields
$$
0\equiv 
(n{-}2)\,\bigl(
\tau_{ij}-\delta_{ij}\sigmabold_0+\psi_{ij}
\bigr)\w\omega_j\w\omega_0
\mod \thetabold,\etabold\,.
$$
By Cartan's Lemma, there must exist functions~$T_{ijk}=T_{ikj}$ so that
$$
\tau_{ij}-\delta_{ij}\sigmabold_0+\psi_{ij}\equiv 2T_{ijk}\,\omega_k
\mod \thetabold,\etabold,\omega_0\,.
\tag20
$$
Symmetrizing~(20) in $i$ and $j$ yields
$$
\tau_{ij}-\delta_{ij}\sigmabold_0 
\equiv \bigl(T_{ijk}+T_{jik}\bigr)\,\omega_k
\mod \thetabold,\etabold,\omega_0\,,
$$
so that there exist functions~$S_{ij}=S_{ji}$ so that
$$
\tau_{ij}-\delta_{ij}\sigmabold_0 
\equiv S_{ij}\,\omega_0 + \bigl(T_{ijk}+T_{jik}\bigr)\,\omega_k
\mod \thetabold,\etabold\,.
$$
Thus, (16) becomes
$$
d\etabold_i\equiv (\sigmabold_i{-}S_{ij}\,\omega_j)\w\omega_0 
+\sigmabold_0\w\omega_i +T_{jik}\,\omega_k\w\omega_j
\mod \thetabold,\etabold\,.
$$
Replacing~$\sigmabold_i$ by $\sigmabold_i{-}S_{ij}\,\omega_j$ will not
affect~(14) (because of the symmetry~$S_{ij}=S_{ji}$) and will, after a
modification of the functions~$F_{ijkl}$ that preserves their symmetries, 
also preserve~(19).  However, this substitution will simplify the above
equation to
$$
d\etabold_i\equiv \sigmabold_i\w\omega_0+\sigmabold_0\w\omega_i 
                                  +T_{jik}\,\omega_k\w\omega_j
\mod \thetabold,\etabold\,.
$$

On the other hand, skewsymmetrizing (20) in $i$ and $j$ yields
$$
\psi_{ij}\equiv \bigl(T_{ijk}-T_{jik}\bigr)\,\omega_k
\mod \thetabold,\etabold,\omega_0\,,
$$
so~(19) becomes
$$
d\etabold_{ij}\equiv \sigmabold_i\w\omega_j-\sigmabold_j\w\omega_i
               +(n{-}2)\,\bigl(T_{ijk}{-}T_{jik}\bigr)\,\omega_k\w\omega_0
                      +{\ts\frac12}F_{ijkl}\,\omega_k\w\omega_l
\mod \thetabold,\etabold\,.
$$

Writing~$S_{ijk} = -S_{ikj} = T_{jik}-T_{kij}$ and using the fact that~
$T_{ijk}=T_{ikj}$, the structure equations for the ideal~$\cI^{(1)}$ 
take the form
$$
\left.
\aligned
d\thetabold_0&\equiv 0\\
d\etabold_0&\equiv (n{-}2)\,\sigmabold_0\w\omega_0+\sigmabold_i\w\omega_i\\
d\etabold_i&\equiv \sigmabold_i\w\omega_0+\sigmabold_0\w\omega_i 
                            -{\ts\frac12}S_{ijk}\,\omega_j\w\omega_k\\
d\etabold_{ij}&\equiv \sigmabold_i\w\omega_j-\sigmabold_j\w\omega_i\\
       &\qquad +(n{-}2)\,S_{kij}\,\omega_k\w\omega_0
                      +{\ts\frac12}F_{ijkl}\,\omega_k\w\omega_l\\
\endaligned\right\}\mod\thetabold,\etabold
\tag21
$$
where the 1@-forms~$\sigmabold_a$ pull back to each fibers of~$Y\to X$ to
be the 1@-forms~$d\sbold_a$.   

Now, equations (21) do not uniquely determine the~$\sigmabold_a$.  For
any functions~$p_a$ on~$Y$, replacing~$\sigmabold_0$ by~$\sigmabold_0
+p_a\,\omega_a$ and~$\sigmabold_i$ by~$\sigmabold_i+(n{-}2)\,p_i\,\omega_0
+p_0\,\omega_i$ will keep the form of~(21) but will modify the functions~
$S_{ijk}$ and $F_{ijkl}$.  In particular, the traced functions~$S_{iij}$
and $F_{ijij}$ will be replaced by~$S_{iij}{+}(n{-}1)\,p_j$ and 
$F_{ijij}{+}2n(n{-}1)p_0$ respectively.  It follows that there is a unique
choice of the~$\sigmabold_a$ so that the torsion functions satisfy the
trace conditions
$$
S_{iij} = 0\qquad\text{and}\qquad F_{ijij}=0.
\tag22
$$
Henceforth, I assume that the~$\sigmabold_a$ have been chosen so that~(21)
and~(22) hold.

Now, I claim that there is at most one integral element of~$\bigl(\cI^{(1)},
\Omega_2\bigr)$ at each point of~$Y$.  In fact, any such integral element~
$E\subset T_yY$ will be defined by equations of the form
$$
\thetabold_i=\etabold_a=\etabold_{ij}=\sigmabold_0-p_a\,\omega_a
=\sigmabold_i-p_{ia}\,\omega_a = 0
$$
 for some unique numbers~$p_a$, and $p_{ia}$.  By~(21), in order that 
$d\etabold_0$ vanish on such an~$E$, it is necessary and sufficient that
$$
p_{i0} = (n{-}2)\,p_i\qquad\text{and}\qquad p_{ij}=p_{ji}\,.
$$
Moreover, again by (21), in order that $d\etabold_i$ also vanish on such 
an~$E$, it is necessary and sufficient that
$$
p_{ij} = \delta_{ij}\,p_0
\qquad\text{and}\qquad 
S_{ijk}(y)=\delta_{ij}\,p_{k}-\delta_{ik}\,p_{j}\,.
$$
However, the trace condition~(22) implies that this last condition can only
hold if
$$
p_i = 0\qquad\text{and}\qquad  S_{ijk}(y)=0.
$$
Thus, the defining relations for the integral element~$E$ take the
 form
$$
\thetabold_i=\etabold_a=\etabold_{ij}=\sigmabold_0-p_0\,\omega_0
=\sigmabold_i-p_{0}\,\omega_i = 0
$$
and these can only define an integral element at points $y\in Y$ where
$S_{ijk}(y)=0$.  Finally, for ~$d\etabold_{ij}$ to vanish
on~$E$, it is necessary and sufficient that~$S_{ijk}(y)=0$ and
$$
 F_{ijkl}(y) - p_0 (\delta_{ik}\delta_{jl}-\delta_{il}\delta_{jk}) = 0.
$$
Again, using the trace condition~(22), this implies $p_0=0$ and
$F_{ijkl}(y)=0$.  

Thus, the subspace~$\Qbold_y\subset T_yY$ defined by
the equations
$$
\thetabold_i=\etabold_a=\etabold_{ij}=\sigmabold_0=\sigmabold_i = 0
\tag23
$$
is the only possible integral element of~$\bigl(\cI^{(1)},\Omega_2\bigr)$
at~$y\in Y$ and this really is an integral element if and only if 
$y$ satisfies
$$
S_{ijk}(y)=F_{ijkl}(y)=0.
\tag24
$$

Since, at each $y\in Y$, the space~$\Qbold_y$ contains the tangents
to the $\Or(n)$@-orbit through~$y$, there is a unique 
$(n{+}1)$@-plane field~$H\subset T\Sigma_2$ that pulls back via the 
projection~$Y\to\Sigma_2$ to be~$\Qbold\subset TY$.  

 Finally, a unit vector field~$\ubold$ on~$U\subset M$ satisfies~(4) 
and~(5) if and only if the bundle~$F_\ubold$ lifts to~$Y$ to be
an integral manifold of~$\bigl(\cI^{(1)},\Omega_2\bigr)$, i.e., so
that it is tangent to the plane field~$\Qbold$.  Since the
image of this lifting under the projection~$Y\to\Sigma_2$ is the
image of the section~$j^2\ubold$ of~$\Sigma_2$, it follows
that $\ubold$ satisfies~(4) and (5) if and only if the image of~$j^2\ubold$ 
is tangent to~$H$. \qed
\enddemo
 
\remark{Remark: Interpretation}
One way of interpreting Theorem~2 is as follows:  All vector fields~$\ubold$ 
that satisfy ~(4) and (5) satisfy a certain total differential equation of 
the form
$$
\nabla^3\ubold = \Rbold\bigl(\ubold,\nabla\ubold,\nabla^2\ubold)
$$
where~$\Rbold$ is a certain nonlinear bundle map from~$J^2(M,SM)$ to 
$TM\otimes T^*M\otimes T^*M\otimes T^*M$.  The formula for~$\Rbold$ in
terms of the metric~$g$ can be written out, but it is unenlightening and, 
moreover, not of much use for doing calculations, so I will not write it 
out here.
\endremark 

\remark{Remark: Integrability}
When $n$ is sufficiently large, for the generic metric~$g$, the equations
$S_{ijk}=F_{ijkl}=0$ on~$Y$ will have no solution in~$Y$.  In this case, 
there will be no harmonic morphism of corank~1 whose domain is an open
subset of~$M$.  

Even when the locus~$Z\subset Y$ defined by~$S_{ijk}=F_{ijkl}=0$ is 
non-empty and smooth, it can well happen that there is no point~$z\in Z$ 
so that $\Qbold_z$ is a subspace of~$T_zZ$.  Again, in this case, there
will be no harmonic morphism of corank~1 whose domain is an open
subset of~$M$. 

Even when there is a smooth submanifold~$Z\subset Y$ that lies in
the locus defined by~$S_{ijk}=F_{ijkl}=0$ with the property that $\Qbold_z$ 
is a subspace of~$T_zZ$ for all~$z\in Z$, it can still happen that~
$\Qbold$ is not an integrable plane field on~$Z$.   

 Finally, when there is a smooth submanifold~$Z\subset Y$ that lies in
the locus defined by~$S_{ijk}=F_{ijkl}=0$ with the property that $\Qbold_z$ 
is a subspace of~$T_zZ$ for all~$z\in Z$ and that~$\Qbold$ is an integrable 
plane field on~$Z$, then~$Z$ will be foliated by integral manifolds of~
$\cI^{(1)}$. Each such integral manifold will correspond to a locally defined
harmonic morphism on some domain in~$M$.
\endremark 

Now, the functions~$S_{ijk}$ and $F_{ijkl}$ on~$Y$ can be expressed explicitly 
in terms of the functions~$\rbold_a$, $\abold_{ij}$, $\sbold_a$, and the
Riemann curvature tensor components~$R_{ijkl}$ and their first covariant 
derivatives on~$F$.  These formulae are very messy in general and it is
difficult to tell much about them.  However, in the case that~$g$ has 
constant sectional curvature, these formulae simplify and constructive
use can be made of them. In the next section, the system~$\cI^{(1)}$
will be examined for such metrics, resulting in a complete classification 
of the harmonic morphisms of corank one with Riemannian domain having 
constant sectional curvature.

\head 3. Harmonic morphisms of corank 1 from space forms  \endhead
In this section, $M$ will always denote a manifold of dimension~$(n{+}1)$ 
that is endowed with a metric~$g$ with constant sectional curvature~$K$.
In this case, the structure equations on the orthonormal frame bundle~
$\pi:F\to M$ simplify to
$$
\align
d\omega_a &= -\omega_{ab}\w\omega_b\,,\\
d\omega_{ab} &= -\omega_{ac}\w\omega_{cb} + K\,\omega_a\w\omega_b\,,\\
\endalign
$$
i.e., $R_{abcd}=K\,\bigl(\delta_{ac}\delta_{bd}-\delta_{ad}\delta_{bc}\bigr)$.

\subhead 3.1. Analysis of the exterior differential system \endsubhead
As in the proof of Theorem~2, define the manifold
$$
X = F\times\bbR\times\bbR^n\times\euso(n)
$$ 
and let~$\rbold_0:X\to\bbR$, $\rbold = (\rbold_i):X\to\bbR^n$, 
and~$\abold = (\abold_{ij}):F\to\euso(n)$ denote the 
projections onto the second third and fourth factors of this product, 
respectively.  As before define 1@-forms~$\thetabold_i$ on~$X$ by 
the formulae
$$
\thetabold_i
  =\omega_{i0}+(n{-}2)\rbold_i\,\omega_0
    +(\rbold_0\,\delta_{ij}-\abold_{ij})\,\omega_j\,.
$$
Next, define 
$$
Y = X\times\bbR\times\bbR^n
$$
and let~$\sbold_0:Y\to\bbR$ and $\sbold = (\sbold_i):Y\to\bbR^n$ be the
projections on the second and third factors respectively.

Define 1-forms~$\etabold_a$ and $\etabold_{ij}$ via the equations~(8)
and~(10) of \S2.2.  Note, however, that the defining equations for~
$\rbold_{ij}$ and~$\bbold_{ijk}$ simplify to 
$$
\align
\bbold_{ijk}&= (n{-}2)\,\bigl(\rbold_i\abold_{jk}+\rbold_j\abold_{ki}
                     +\rbold_k\abold_{ji}\bigr)\\
-(n{-}2)\,\rbold_{ij} 
&= \delta_{ij}\,\bigl({\rbold_0}^2+K\bigr)
      +\abold_{ik}\abold_{kj} +(n{-}2)^2\,\rbold_i\rbold_j\,.\\
\endalign   
$$

\subsubhead 3.1.1. First torsion equations \endsubsubhead
With great effort, 1@-forms~$\sigmabold_a$ and functions~$S_{ijk}$
and~$F_{ijkl}$ can be explicitly computed on~$Y$ so that equations~(21)
and~(22) hold. In particular, this computation yields
$$
(n{-}2)\,S_{ijk} 
= (n{-}1)\bigl(\pbold_j\,\abold_{ki}-\pbold_k\,\abold_{ji}
                       -2\pbold_i\,\abold_{jk}\bigr)
   -3\bigl(\delta_{ij}\,\pbold_l\,\abold_{kl}
         -\delta_{ik}\,\pbold_l\,\abold_{jl}\bigr)
\tag$1_{ijk}$
$$
where
$$
\pbold_i = \sbold_i - (n{-}2)\,\rbold_0\rbold_i\,.
\tag$2_{i}$
$$
The formula for~$F_{ijkl}$ is much more complicated.  Thankfully, it
will not be needed.

Let~$Z\subset Y$ be the locus defined by the equations~$S_{ijk}=0$.  It
is clear from equation~(1) that $Z$ contains the two loci~$Z_1$ and
$Z_2$ where
$$
\align
Z_1 &= \bigl\{z\in Y \mid \abold_{ij}(z)=0\ \text{for all~$i,j$\ }\bigr\},\\
Z_2 &= \bigl\{z\in Y \mid \pbold_{i }(z)=0\ \text{for all~$i$\ }\bigr\}.\\
\endalign
$$
I claim that, in fact, $Z = Z_1\cup Z_2$.  To prove this, it is enough to 
show, for any $p=(p_i)\in\bbR^n$ and any skew-symmetric matrix~$a=(a_{ij})$, 
that the equations
$$
0= (n{-}1)\bigl(p_j\,a_{ki}-p_k\,a_{ji}-2\,p_i\,a_{jk}\bigr)
   -3\bigl(\delta_{ij}\,p_l\,a_{kl}-\delta_{ik}\,p_l\,a_{jl}\bigr)
\tag$1'_{ijk}$
$$
hold for all~$i$, $j$, and $k$ only if either $p_i=0$ for all~$i$
or else~$a_{ij}=0$ for all $i$ and $j$.  To show this, set~$q_i=a_{ij}\,p_j$.
Note that, because of the skew-symmetry of~$a$, it follows that
$q_i\,p_i = a_{ij}\,p_jp_i = 0$.  Multiplying~$(1'_{ijk})$ by~$p_k$ 
and summing over~$k$ yields
$$
0= (n{-}1)\bigl(-p_j\,q_i-(p_k)^2\,a_{ji}-2\,p_i\,q_j\bigr) +3\,q_j\,p_i
$$
which can be rearranged to give
$$
(n{-}1)\,(p_k)^2\,a_{ji} = -(n{-}1)\,p_j\,q_i - (2n{-}5)\,p_i\,q_j
\tag$3_{ij}$
$$
Symmetrizing~$(3_{ij})$ in~$i$ and~$j$ yields the relation~$0=p_iq_j+p_jq_i$, 
valid for all~$i$ and~$j$. Multiplying this last relation by ~$p_i$ and
summing over~$i$ yields~$(p_i)^2\,q_j = 0$.  Thus, unless~$p_i=0$ for all~
$i$, then ~$q_j = 0 $ for all~$j$, and, thence, by the equations~$(3_{ij})$, 
it would follow that~$a_{ij}=0$ for all ~$i$ and $j$.  
Thus, the claim is proved.

Now, from the proof of Theorem~2, any integral manifold of~$\bigl(\cI^{(1)},
\Omega_2\bigr)$ in~$Y$ must lie in the locus~$Z = Z_1\cup Z_2$.
Since all of the connected integral manifolds of the system~$\cI^{(1)}$ are 
derived from solutions of the harmonic morphism equations and hence must
be real analytic, it follows that any connected integral manifold 
of~$\cI^{(1)}$ must lie in either~$Z_1$ or $Z_2$ (or possibly, both).  
Thus, there are two types of solutions to be studied:  The first type
are the integral manifolds that lie in~$Z_1$ and the second type are
the integral manifolds that lie in~$Z_2$.

\subsubhead 3.1.2. Integral manifolds of the first type \endsubsubhead
Since~$Z_1$ is defined by the equations~$\abold_{ij}=0$, it follows from the 
definitions and the simplified formula for~$\bbold_{ijk}$ that, on~$Z_1$, 
the formula for~$\etabold_{ij}$ simplifies to
$$
\etabold_{ij} = \sbold_i\,\omega_j - \sbold_j\,\omega_i\,.
$$
In particular, any integral manifold of~$\bigl(\cI^{(1)},\Omega_2\bigr)$ 
in~$Z_1$ must lie in the sublocus~$Z_{11}\subset Z_1$ defined by the 
equations~$\abold_{ij}=\sbold_{j}=0$.  

Now, on~$Z_{11}$, the ideal~$\cI^{(1)}$ is generated by the 1@-forms
$$
\align
\thetabold_i&=\omega_{i0}+(n{-}2)\rbold_i\,\omega_0+\rbold_0\,\omega_i\,,\\
  \etabold_0&=d\rbold_0+ (n{-}2)\,\sbold_0\,\omega_0\,,\\
  \etabold_i&=d\rbold_i - \rbold_j\,\omega_{ji} 
                 - (n{-}1)\,\rbold_0\rbold_i\,\omega_0
                 +(n{-}2)\,\rbold_i\rbold_j\,\omega_j
       +\bigl(\sbold_0+{\ts\frac1{n-2}}({\rbold_0}^2+K)\bigr)\,\omega_i
\endalign
$$
and their exterior derivatives.  With a slight notation change, 
this becomes
$$
\aligned
\thetabold_i&=\omega_{i0}+(n{-}2)\rbold_i\,\omega_0+\rbold_0\,\omega_i\,,\\
  \etabold_0&=d\rbold_0
               +\bigl((n{-}2)\,\sbold_0-{\rbold_0}^2-K\bigr)\,\omega_0\,,\\
  \etabold_i&=d\rbold_i - \rbold_j\,\omega_{ji} 
                 - (n{-}1)\,\rbold_0\rbold_i\,\omega_0
                 +(n{-}2)\,\rbold_i\rbold_j\,\omega_j
                 +\sbold_0\,\omega_i\,.
\endaligned
\tag4
$$
A short computation yields
$$
\left.
\aligned
d\thetabold_i &\equiv   0\\
  d\etabold_0 &\equiv  \sigmabold_0\w\omega_0\\
  d\etabold_i &\equiv  \sigmabold_0\w\omega_i
               -2(n{-}1)(n{-}2)\,\rbold_0\rbold_i\rbold_k\,\omega_k\w\omega_0
\endaligned
\right\} \mod\thetabold,\etabold,
\tag5
$$
where
$$
\sigmabold_0 
= d\sbold_0 
    -\rbold_0\bigl(n\,\sbold_0+(n{-}2){\rbold_k}^2\bigr)\,\omega_0
    +\bigl((n{-}2)\sbold_0-{\rbold_0}^2-K\bigr)\,\rbold_k\,\omega_k          
$$
 From~(5) and the real analyticity mentioned before, it follows that 
any connected integral manifold of~$\bigl(\cI^{(1)},\Omega_2\bigr)$ 
that lies in~$Z_{11}$ must lie in either the sublocus~$Z_{111}\subset Z_{11}$ 
defined by the equation~$\rbold_0= 0$ or the sublocus~$Z_{112}\subset Z_{11}$ 
defined by the equations~$\rbold_i=0$ for all~$i=1,\dots n$.  

Now, pulled back to~$Z_{111}$, the 1@-form~$\etabold_0$ simplifies to
$\etabold_0 = \bigl((n{-}2)\,\sbold_0-K\bigr)\,\omega_0$, so any 
integral manifold of~$\bigl(\cI^{(1)},\Omega_2\bigr)$ that lies in~$Z_{111}$ 
must lie in the sublocus~$Z_{1111}\subset Z_{111}$ defined by the equation
$\sbold_0 = K/(n-2)$.  Now, on~$Z_{1111}$, the generators of~$\cI^{(1)}$
simplify to
$$
\aligned
\thetabold_i&=\omega_{i0}+(n{-}2)\rbold_i\,\omega_0\,,\\
  \etabold_i&=d\rbold_i - \rbold_j\,\omega_{ji}
                 +(n{-}2)\,\rbold_i\rbold_j\,\omega_j
                 +{\ts\frac1{(n-2)}}\,K\,\omega_i\,.
\endaligned
\tag6
$$
and these have structure equations
$$
d\thetabold_i\equiv d\etabold_i\equiv 0 \mod \thetabold,\etabold.
\tag7
$$
Thus, the system~$\cI^{(1)}$ on~$Z_{1111}$ is a Frobenius Pfaffian system,
so that~$Z_{1111}$ is foliated by integral manifolds of~
$\bigl(\cI^{(1)},\Omega_2\bigr)$.  I will call the foliations of~$M$
corresponding to these integral manifolds, foliations of Type~0.  
They will be analyzed below.

Now examine the system~$\bigl(\cI^{(1)},\Omega_2\bigr)$ on the
sublocus~$Z_{112}\subset Z_{11}$.  On this sub\-locus, the 1@-forms~
$\etabold_i$ simplify to~$\etabold_i = \sbold_0\,\omega_i$, so all of the 
integral manifolds of~$\bigl(\cI^{(1)},\Omega_2\bigr)$ in~$Z_{112}$ must, 
in fact, lie in the sublocus~$Z_{1121}\subset Z_{112}$ defined by the 
equation~$\sbold_0=0$.  On~$Z_{1121}$, the system~$\cI^{(1)}$ has generators
$$
\aligned
\thetabold_i&=\omega_{i0}+\rbold_0\,\omega_i\,,\\
  \etabold_0&=d\rbold_0-\bigl({\rbold_0}^2+K\bigr)\,\omega_0\,,\\
\endaligned
\tag8
$$
and these satisfy the structure equations
$$
d\thetabold_i\equiv d\etabold_0\equiv 0 \mod \thetabold,\etabold_0\,.
$$
In particular, the system~$\cI^{(1)}$ is a Frobenius Pfaffian system
on~$Z_{1121}$ and so it is foliated by integral manifolds of~
$\bigl(\cI^{(1)},\Omega_2\bigr)$.  I will call the foliations of~$M$
corresponding to these integral manifolds, foliations of Type~1.  
They will be analyzed below.

\subsubhead 3.1.3. Integral manifolds of the second type \endsubsubhead
Now consider the locus~$Z_2\subset Z$ defined by the equations
$\sbold_i = (n{-}2)\,\rbold_0\rbold_i$. Calculation now shows that
on~$Z_2$, the relation
$$
d\etabold_0
\equiv (n{-}2)\,\sigmabold_0\w\omega_0
    - (n-2)\,(\sbold_0-{\rbold_0}^2)\abold_{ij}\,\omega_i\w\omega_j
\mod \thetabold,\etabold.
$$
holds, where $\sigmabold_0\equiv d\sbold_0\mod\omega$.  It follows that
any integral manifold of~$\bigl(\cI^{(1)},\Omega_2\bigr)$ in~$Z_{2}$
must lie in either the sublocus~$Z_{21}\subset Z_2$ defined by the equation~
$\sbold-{\rbold_0}^2=0$ or else in the sublocus~$Z_1\cap Z_2\subset Z_2$
defined by the equations~$\abold_{ij}=0$.  Since all of the integral manifolds
of~$\bigl(\cI^{(1)},\Omega_2\bigr)$ that lie in~$Z_1$ have already been
 found, this second case will be discarded.

Thus, consider the system~$\bigl(\cI^{(1)},\Omega_2\bigr)$ on~$Z_{21}$. 
The 1@-form generators of this system are
$$
\aligned
\thetabold_i
  &=\omega_{i0}+(n{-}2)\rbold_i\,\omega_0
    +(\rbold_0\,\delta_{ij}-\abold_{ij})\,\omega_j\,,\\
\etabold_0 &= d\rbold_0+(n{-}2)\,\rbold_0
             \bigl(\rbold_0\,\omega_0+\rbold_i\,\omega_i\bigr)\,,\\
\etabold_i &= d\rbold_i - \rbold_j\,\omega_{ji}
               -\bigl(\rbold_0\rbold_i
                          -\abold_{ij}\rbold_j\bigr)\,\omega_0 
               +(\rbold_0\,\abold_{ij}-\rbold_{ij}
                   +{\rbold_0}^2\,\delta_{ij})\,\omega_j\,,\\
\etabold_{ij} &= d\abold_{ij}
                   - \abold_{kj}\,\omega_{ki}-\abold_{ik}\,\omega_{kj}\\
          &\qquad - n\,\rbold_0\abold_{ij}\,\omega_0
                  - (n{-}2)\,\bigl(\bbold_{ijk}-\rbold_0\rbold_i\,\delta_{jk}
                       +\rbold_0\rbold_j\,\delta_{ik}\bigr)\,\omega_k\,,\\
\endaligned
\tag9
$$
where, for simplicity, $\bbold_{ijk}$ has been redefined to be
$$
\bbold_{ijk}=\rbold_i\abold_{jk}+\rbold_j\abold_{ki}+\rbold_k\abold_{ji}\,.
$$

On~$Z_{21}$, the congruences~$d\thetabold_i\equiv d\etabold_0\equiv0
\mod \thetabold,\etabold$ hold, but
$$
d\etabold_i
   \equiv 2(n{-}1){\rbold_0}^2\,\omega_i\w\bigl(\rbold_k\,\omega_k\bigr)
    \mod \omega_0,\thetabold,\etabold.
$$
It follows that any connected integral manifold of~$\bigl(\cI^{(1)},
\Omega_2\bigr)$ that lies in~$Z_{21}$ must lie in either the sublocus
$Z_{211}\subset Z_{21}$ defined by the equation~$\rbold_0=0$ or the
sublocus~$Z_{212}\subset Z_{21}$ defined by the equations~$\rbold_i=0$.  

On $Z_{211}$, the generators of~$\cI^{(1)}$ simplify to 
$$
\aligned
\thetabold_i 
   &=\omega_{i0}+(n{-}2)\rbold_i\,\omega_0-\abold_{ij}\,\omega_j\,,\\
\etabold_i &= d\rbold_i - \rbold_j\,\omega_{ji}
                        + \abold_{ij}\rbold_j\,\omega_0 
                        - \rbold_{ij}\,\omega_j\,,\\
\etabold_{ij} &= d\abold_{ij}
                   - \abold_{kj}\,\omega_{ki}-\abold_{ik}\,\omega_{kj}
  - (n{-}2)\,(\rbold_i\abold_{jk}{+}\rbold_j\abold_{ki}{+}\rbold_k\abold_{ji})
            \,\omega_k\,,\\
\endaligned
\tag10
$$
where
$$
-(n{-}2)\,\rbold_{ij} = K\delta_{ij}+\abold_{ik}\abold_{kj}
                        + (n{-}2)^2\,\rbold_i\rbold_j\,.
$$
Now computation yields
$$
d\thetabold_i\equiv d\etabold_i\equiv d\etabold_{ij}\equiv 0
\mod\thetabold,\etabold.
$$
Thus, on~$Z_{211}$, the system~$\cI^{(1)}$ is a Frobenius Pfaffian system.
It follows that~$Z_{211}$ is foliated by integral manifolds of~
$\bigl(\cI^{(1)},\Omega_2\bigr)$.  I will call the foliations of~$M$
corresponding to these integral manifolds, foliations of Type~2.  Note
that~$Z_{1111}$ is a submanifold of~$Z_{211}$ and that the integral
manifolds of Type~0 are special cases of those of Type~2.
They will be analyzed below.

 Finally, consider the system~$\cI^{(1)}$ on the locus~$Z_{212}$.  On
this submanifold, the formula for~$\etabold_i$ simplifies to
$$
\etabold_i = (\rbold_0\,\abold_{ij}-\rbold_{ij}
                   +{\rbold_0}^2\,\delta_{ij})\,\omega_j
$$
Since~$\rbold_{ij}=\rbold_{ji}$, the vanishing of~$\etabold_i$ on a connected
integral manifold of~$\bigl(\cI^{(1)},\Omega_2\bigr)$ that lies in~$Z_{212}$ 
implies that the functions~$\rbold_0\,\abold_{ij}$ must also vanish on such 
an integral manifold.  Thus, such an integral manifold must either satisfy~
$\rbold_0=0$, in which case, it lies in~$Z_{211}$ and so has already been 
accounted for, or else satisfy~$\abold_{ij}=0$, in which case, it lies
in~$Z_1$ and so has already been accounted for.

\subsubhead 3.1.4. Conclusions \endsubsubhead
In summary, there are two types of integral manifolds of~$\bigl(\cI^{(1)},
\Omega_2\bigr)$:  When~$\cI^{(1)}$ is pulled back to either $Z_{1121}$ or~
$Z_{211}$, it becomes a Fro\-benius system whose leaves are maximal integral 
manifolds of~$\bigl(\cI^{(1)},\Omega_2\bigr)$.  Moreover every connected 
integral manifold of~$\bigl(\cI^{(1)},\Omega_2\bigr)$ is an open subset of
one of these leaves.

\subhead 3.2. The local classification \endsubhead
 Finally, the main result can be stated and proved.  Note that this
theorem does not require any assumption of compactness or completeness.

\proclaim{Theorem 3}
Let~$\bigl(M^{n+1},g\bigr)$ be a $1$@-connected manifold of constant sectional 
curvature~$K$ and suppose that~$\phi:\bigl(M^{n+1},g\bigr)\to
\bigl(N^{n},h\bigr)$ be a submersive harmonic morphism with connected fibers.  
Then either there exists a Killing field~$\Xbold$ tangent to the fibers
of~$\pi$ or else the fibers of~$\pi$ are geodesics and~$M$ is foliated by 
totally umbilic hypersurfaces that are orthogonal to the fibers of~$\pi$.
\endproclaim

\demo{Proof}
Since~~$M$ is $1$@-connected and the fibers of~$\pi$ are connected, it
 follows that the fibers are orientable, i.e., that there exists a unit
vector field~$\ubold$ on~$M$ whose fibers are the integral curves of~$\ubold$
in~$M$.  Using the notation established in the proof of Theorem~2, let 
$F_\ubold\subset F$ be the subbundle of the $g$@-orthonormal frame bundle
of~$M$ consisting of those frames~$(m;u)$ so that~$u_0 = \ubold(m)$.

Then there exist functions~$r_0$, $r_i$, and $a_{ij}=-a_{ji}$ on~$F_\ubold$
so that 
$$
\omega_{i0} = -(n{-}2) r_i\,\omega_0-(r_0\,\delta_{ij}-a_{ij})\,\omega_j\,.
$$
Moreover, these functions satisfy
$$
d\bigl( r_0\,\omega_0 + r_i\,\omega_i \bigr) = 0.
$$

By the analysis in \S3.1, there are two possibilities.  Either $r_i=a_{ij}=0$
and $dr_0 = ({r_0}^2 + K)\,\omega_0$ (if the corresponding integral manifold
is of Type~1) or else $r_0=0$ and the functions~$r_i$ and $a_{ij}$ satisfy
the equations
$$
\aligned
dr_i &=  r_j\,\omega_{ji}-a_{ij}r_j\,\omega_0 + r_{ij}\,\omega_j\,,\\
da_{ij} &=  a_{kj}\,\omega_{ki} + a_{ik}\,\omega_{kj}
            + (n{-}2)\,(r_ia_{jk}{+}r_ja_{ki}{+}r_ka_{ji})\,\omega_k\,,\\
\endaligned
$$
where
$$
-(n{-}2)\,r_{ij} = K\delta_{ij}+a_{ik}a_{kj}
                        + (n{-}2)^2\,r_ir_j\,,
$$
(if the corresponding integral manifold is of Type~2).

Suppose that the first possibility holds.  Then, since~$r_i=0$, it follows
 from the discussion in~\S1.1.3 that the fibers of~$\phi$ are geodesics (since
they have vanishing mean curvature vector).  Moreover, now the identity
$\omega_{i0} = -r_0\,\omega_i$ holds, so it follows that ~$d\omega_0=0$, 
i.e., the leaves of~$\omega_0=0$ are the frame bundles of hypersurfaces in~$M$
that are orthogonal to the fibers of~$\phi$.  By the structure equations, 
the second fundamental form~$H$ of each such leaf is the restriction to
that leaf of the tensor
$$
\Hbold = \omega_{i0} \circ \omega_i = -r_0\,\omega_i\circ\omega_i = -r_0\,g'
$$
where~$g'$ is the usual orthogonal projection of the metric~$g$.  Since~
$\Hbold$ is a scalar multiple of~$g'$, the induced metric on these orthogonal 
hypersurface leaves, it follows that each such hypersurface is totally
umbilic.

Suppose now that the second possibility holds.  Since~$r_0=0$, it follows
that~$\rho = \rho' = r_i\,\omega_i$ is a closed 1-form that is well-defined
on~$M$.  Since~$M$ is connected and simply connected, there exists a smooth
positive function~$r$ on~$M$, unique up to a constant scalar multiple, so 
that~$r^{-1}\,dr = \rho$.  Set $\Xbold = r^{n-2}\,\ubold$.   Now a calculation
using the equation~$\omega_{i0} = -(n{-}2) r_i\,\omega_0+a_{ij}\,\omega_j$  
yields
$$
\Lie_\Xbold g = 0.
$$
Thus~$\Xbold$ is a Killing vector field for~$g$, as desired. \qed
\enddemo

\remark{Remark}  The two possibilities are not quite mutually exclusive.
There is essentially one local example that falls under both Types.  This is
when~$r_0 = r_i = a_{ij} = 0$.  By the structure equation~$dr_0 = 
({r_0}^2 + K)\,\omega_0$, this can only happen when~$K=0$.  Then the foliation
of~$\phi$@-fibers is by parallel geodesics and the orthogonal hypersurfaces
are totally geodesic.  In this case, the morphism~$\phi$ is locally equivalent
to the standard linear orthogonal projection~$\phi:\bbR^{n+1}\to\bbR^n$.
\endremark

\subhead 3.3. Examples and further results \endsubhead
Theorem~3 forms the basis of a classification of the harmonic morphisms
of corank~1 for which the domain has constant sectional curvature.   

\subsubhead 3.3.1.  Umbilic morphisms \endsubsubhead
Because of the nature of the orthogonal foliation, I will say that
a harmonic morphism~$\phi:\bigl(M^{n+1},g\bigr)\to\bigl(N^n,h\bigr)$ whose
corresponding integral manifold is of Type~1 is an {\it umbilic morphism}.
Thus, in this case, the fibers of~$\phi$ are the geodesics orthogonal
to a foliation of~$M$ by parallel totally umbilic hypersurfaces.  To specify
the map~$\phi$ up to obvious equivalences, it is evidently enough to
specify a single totally umbilic hypersurface in~$M$.  The number of
possible cases depends on the sign of the sectional curvature~$K$.  For
simplicity, I will only treat the cases~$K = 1,0,-1$.

If~$K=1$ and~$M=S^{n+1}$ with its standard metric, then every totally 
umbilic hypersurface is an $n$@-sphere of radius at most~1.  The family
of parallel hypersurfaces consist of the parallel hyperspheres, one of
which is a great~$n$-sphere.  Using the standard inclusions of~$S^p$
into~$\bbR^{n+1}$, the formula for the map~$\phi$ becomes
$$
\phi\bigl(x_0,\dots,x_{n+1}\bigr)
=\frac{\bigl(x_0,\dots,x_{n}\bigr)}{\sqrt{1-{x_{n+1}}^2}}\,,
$$
which is undefined at the `poles', $x_{n+1} = \pm1$.  

If~$K=0$ and $M=\bbR^{n+1}$, there are two types of umbilic foliations.  
The first type is by parallel planes, leading to the harmonic morphism
$$
\phi\bigl(x_0,\dots,x_{n+1}\bigr)=\bigl(x_0,\dots,x_{n}\bigr)\,,
$$
and the second is the radial projection~$\phi:\bbR^{n+1}\setminus{\bold 0}
\to S^n$ given by 
$$
\phi\bigl(x_1,\dots,x_{n+1}\bigr)
=\frac{\bigl(x_1,\dots,x_{n}\bigr)}{\sqrt{{x_1}^2+\dots+{x_{n+1}}^2}}\,.
$$

If~$K=-1$, and~$M=D^{n+1}\subset\bbR^{n+1}$, the hyperbolic ball,
then there are three types of umbilic foliations, leading to three types of 
umbilic morphisms.  The first type is when one of the umbilic hypersurfaces 
is totally geodesic.  Then the parallel hypersurfaces are all totally umbilic 
with principal curvatures less than~1 in absolute value. 
The parallel mapping identifying any two such
hypersurfaces is a homothety and the quotient metric is of constant
negative sectional curvature.  Thus, the corresponding~$\phi$ is a harmonic
morphism from~$D^{n+1}$ to~$D^n$.  The second type is when all of the
parallel umbilic hypersurfaces have principal curvature equal to $+1$ or $-1$.
Of course, in this case, the hypersurfaces are all horocycles tangent at
a unique common point on the ideal boundary of~$D^{n+1}$.  The metric on such
a horocycle is the flat metric, so this gives rise to a harmonic morphism
$\phi:D^{n+1}\to\bbR^n$.  (This example is more simply seen in coordinates
as the linear projection from the upper half space model of hyperbolic space 
to the boundary plane.)  Finally, the third type is when all of the
parallel umbilic hypersurfaces have principal curvatures of absolute
value greater than~$1$.  In this case, the hypersurfaces are the level sets
of the distance function from some fixed point~$z\in D^{n+1}$ and hence
are isometric to standard spheres~$S^n$ of varying radii.  The corresponding
harmonic morphism is well-defined as a map~$\phi:D\setminus\{z\}\to S^n$.  

By passing to quotients and covers and so forth, every example with~$M$
complete or compact can be constructed from these examples.  Note that by 
taking appropriate quotients by discrete subgroups, one can construct compact 
examples only by starting with the harmonic morphisms from~$\bbR^{n+1}$
or~$D^{n+1}$ to~$\bbR^n$ since these are the only ones for which~$r_0$
can satisfy the equation~$dr_0 = ({r_0}^2+K)\,\omega_0$ on a compact fiber
of~$\phi$.

All of these examples of umbilic morphisms are due to Gudmundsson,
who also proved a characterization theorem~\cite{Gu1, Theorem~3.6} asserting
that a non-constant harmonic morphism~$\phi:\bigl(M^{n+1},g\bigr)\to
\bigl(N,h\bigr)$ where both~$g$ and~$h$ have constant sectional curvature 
that has the additional property that the conformal factor~$r$ is constant
on curves in~$M$ perpendicular to the fibers of~$\phi$ must be one of
these examples or else~$r$ is actually constant.

\subsubhead 3.3.2.  Isometric quotient morphisms \endsubsubhead
 Finally, consider the case where there exists a non-zero Killing field~
$\Xbold$ on~$M$ whose integral curves are the fibers of the harmonic 
morphism~$\phi$. From the discussion of this example in~\S2.1, the metric~$h$ 
on the target manifold~$N$ must satisfy~
$\phi^*(h)=c\,|\Xbold|^{2/(n-2)}\,g'$ for some constant~$c>0$, where, 
as usual, $g'$ denotes the part of the metric~
$g$ that is orthogonal to the fibers of~$\phi$.  If $M$ is oriented and~
${*}_g1$ denotes the volume form of~$g$, then a computation shows that the 
pulled-back volume form for~$h$ must have the form
$$
\phi^*({*}_h1) = c^{n/2}\,|\Xbold|^{2/(n-2)} \bigl(\Xbold\lefthook{*}_g1\bigr).
$$ 
Note that the $n$@-form on the right hand side of this equation is 
smooth away from the zero locus~$Z\subset M$ of~$\Xbold$.   

Because~$\Xbold$ is a non-zero Killing field, $Z$ must be a smooth, 
proper submanifold of even codimension, say,~$2q\le n{+}1$.   In fact, in a 
neighborhood of any~$z\in Z$, there exists a geodesic normal coordinate 
system $y_1,\dots,y_{n+1}$ on a $z$@-neighborhood~$U$ so that~$Z\cap U$ is 
defined by the equations~$y_1 = \dots = y_{2q} = 0$ and~$\Xbold$ has the
 form
$$
\Xbold =  m_1\left(y_1\,\frac{\partial}{\partial y_2}
                  -y_2\,\frac{\partial}{\partial y_1}\right)
        + \dots
        + m_q\left(y_{2q-1}\,\frac{\partial}{\partial y_{2q}}
                  -y_{2q}\,\frac{\partial}{\partial y_{2q-1}}\right)
$$
where~$m_i>0$.  It follows that the expression~
$|\Xbold|^{2/(n-2)} \bigl(\Xbold\lefthook{*}_g1\bigr)$ cannot be smooth
on a neighborhood of~$z$ unless~$n=3$.  Of course, when~$n=3$, this expression
actually is smooth, a fact that will be used below.  Moreover, when~$n=3$,
the formula for the `push-down' metric simplifies to
$$
\phi^*(h) = c\,\left(\,|\Xbold|^2 g - (\Xbold^\flat)^2\,\right),
$$
where~$\Xbold^\flat$ represents the $g$@-dual $1$@-form to~$\Xbold$.

When~$M^{n+1}$ is a simply connected manifold and~$g$ is complete
with constant sectional curvature~$K$, there are many Killing vector fields,
each leading to a different harmonic quotient. Of course, these quotients
will usually have singularities, since the space of integral curves of
a general Killing vector field on~$M$ will not, in general, carry the
structure of a smooth manifold.  Thus, it will usually be necessary to 
restrict to an open set~$U\subset M$ on which the integral curves of~$\Xbold$
 form an amenable foliation before one can construct the target manifold.

\example{Example 2}  {\it Spherical quotients.}
Consider the case~$K=1$, so that~$M=S^{n+1}\subset\bbR^{n+2}$ 
(embedded as the standard unit sphere~${x_0}^2+\dots+{x_{n+1}}^2 = 1$).  
Then, in order for the integral curves of~$\Xbold$ to be closed (so that 
the space of integral curves endowed with the quotient topology is at least 
Hausdorff), it must be conjugate in the rotation group to a multiple of a 
vector field of the form
$$
\Xbold =  m_0\left(x_0\,\frac{\partial}{\partial x_1}
                  -x_1\,\frac{\partial}{\partial x_0}\right)
        + \dots
        + m_k\left(x_{2k}\,\frac{\partial}{\partial x_{2k+1}}
                  -x_{2k+1}\,\frac{\partial}{\partial x_{2k}}\right)
$$
 for some integer $k$ satisfying~$2k\le n$ and some positive integers~
$m_0\le\dots\le m_k$ with greatest common divisor equal to 1.   Take
$\Xbold$ to have this form.

The generic integral curve of~$\Xbold$ is of period~$2\pi$ and there will be 
ramified integral curves (that are not fixed points) unless~$m_i=1$ for 
all~$i$.  The fixed point set~$Z\subset S^{n+1}$ will be empty unless~
$2k<n$, in which case~$Z = S^{n-2k-1}$ is defined by the equations~
$x_0 = x_1 = \dots = x_{2k+1}=0$.  Note that, away from~$Z$, the 
$\Xbold$@-orthogonal part~$g'$ of the induced metric~$g$ is smooth, as is
the tensor~$|\Xbold|^{2/(n-2)}\,g'$.  

If~$R\subset S^{n+1}$ denotes the
closed subset of~$S^{n+1}$ consisting of the ramified integral curves of~
$\Xbold$ plus the fixed points, then setting~$N = (S^{n+1}\setminus R)/\Xbold$,
with quotient projection~$\phi: S^{n+1}\setminus R\to N$,
defines a smooth $n$@-manifold on which there exists a unique Riemannian
metric~$h$ that satisfies~$\phi^*(h) = |\Xbold|^{2/(n-2)}\,g'$.  With
respect to this metric,~$\phi:\bigl(S^{n+1}\setminus Z, g\bigr)
\to\bigl(N,h\bigr)$ is a harmonic morphism.

Note that since~$|\Xbold|^{2/(n-2)}\,g'$ does not vanish on~$R\setminus Z$,
it is not possible to extend~$N$ as a smooth manifold in such a way that 
$\phi$ can be extended over any part of~$R\setminus Z$ as a smooth mapping.
Thus, if~$\phi$ is to extend globally to~$S^{n+1}$ as a harmonic morphism
to a smooth $n$@-manifold, then~$R=Z$, i.e., $m_i=1$ for all~$i$.  
Assume this from now on, so that~$R=Z$.

Now, if $Z$ is empty, i.e., if $n=2k$, then~$\Xbold$ generates the
standard circle action on~$S^{2k+1}$ and the quotient manifold is~
$N=\bbC\bbP^k$, with the Fubini-Study metric.

On the other hand, if~$Z=S^{n-2k-1}$ is non-empty, it consists of either 
1 or 2 components.  Since the expression~$|\Xbold|^{2/(n-2)}\,\bigl(\Xbold
\lefthook {*}_g1\bigr)$ is not smooth along~$Z$ unless $n=3$, it follows that 
there are at most two cases to consider:  The cases~$(n,k)=(3,0)$ and 
$(n,k)=(3,1)$.  In the first case, $N$ can be completed to a topological
manifold~$\bar N$ by adding a single point corresponding to the fixed
point set~$Z = S^2$ and, in the second case, $N$ can be completed to a 
topological manifold~$\bar N$ by adding two points, corresponding to the two 
points of~$Z = S^0$.  However, in neither case can the smooth structure on~$N$
be extended across the extra point(s) so that the metric~$h$ extends as
a smooth metric on~$\bar N$.  In fact, in each case, the metric on~$N$ is
invariant under an action of $\SO(3)$, so that in a neighborhood of one of the
extra points, `polar coordinates' can be introduced so that~$s$ represents
the distance from the singular point and a neighborhood of the extra point
can be written in the form~$(0,\epsilon)\times S^2$ with the metric taking
the form
$$
ds^2 + f(s)\,\sigma
$$
where $\sigma$ is the standard metric of curvature~1 on~$S^2$.  In the
case that~$k=0$, the function~$f$ has the form~$f(s) = s(s+1)^2(s+2)$
near~$s=0$ (the extra point) while in the case that $k=1$, the function~$f$ 
has the form~$f(s) = s^2(1-s/2)^2$ near $s=0$.  In neither case is $f$ the
square of a smooth, odd function of~$s$ whose derivative at $s=0$ equals~$1$,
a necessary and sufficient condition for such a metric to represent 
a smooth, rotationally invariant metric in polar coordinates.

The conclusion of all this discussion is that the only non-constant
harmonic morphism whose domain is the entire~$S^{n+1}$ with its standard
metric and whose range is a smooth $n$@-manifold is the standard fibering
$\phi:S^{2k+1}\to\bbC\bbP^k$.  Note that this analysis gives a more
comprehensive solution to the problem discussed in~\cite{Do}.
\endexample

\remark{Remark: Other examples}  Using the above discussion as a guide, it
is possible to write down an $\SO(4)$@-invariant metric~$g$ on~
$S^4\subset\bbR^5$ that admits a Killing field~$\Xbold$ with two 
isolated fixed points so that the quotient space~$N$ carries the structure 
of a smooth manifold diffeomorphic to~$S^3$ endowed with a smooth metric~$h$, 
so that the leaf projection~$\phi:\bigl(S^4,g\bigr)\to\bigl(S^3,h\bigr)$ is a
harmonic morphism with exactly two singular points.  However, the metric~$g$
cannot have constant sectional curvature.
\endremark

\example{Example 3}  {\it Euclidean quotients.}
The  non-zero Killing vector fields on~$\bbR^{n+1}$ can be divided 
into two types, those without fixed points and those with fixed points.

If~$\Xbold$ is a Killing field on~$\bbR^{n+1}$ without fixed points, then,
up to a constant multiple, it it conjugate by an Euclidean motion to the
vector field
$$
\Xbold = \frac{\partial}{\partial x_0}
         + m_1\left(x_1\,\frac{\partial}{\partial x_2}
                   -x_2\,\frac{\partial}{\partial x_1}\right)
         + \dots
         + m_q\left(x_{2q-1}\,\frac{\partial}{\partial x_{2q}}
                   -x_{2q}\,\frac{\partial}{\partial x_{2q-1}}\right)
$$
 for some positive real numbers~$m_1,\dots,m_q$ where~$q\le n/2$.  From
now on, assume that~$\Xbold$ has this form.  Note that all of the integral
curves of~$\Xbold$ intersect the hyperplane~$x_0=0$ exactly once, and do so
transversely, so that this hyperplane can be taken as the model space for
the quotient manifold~$N$.  It is not difficult to see that the metric~$h$ 
on~$N$ (canonically determined up to a constant scalar factor) that makes 
the leaf projection~$\phi:\bbR^{n+1}\to N=\bbR^n$ a harmonic morphism  
is never of constant sectional curvature unless~$q=0$, in which case~$\phi$
is simply the linear orthogonal projection~$\phi:\bbR^{n+1}\to\bbR^n$.

If~$\Xbold$ is a Killing field on~$\bbR^{n+1}$ with fixed points, then,
it is conjugate by an Euclidean motion to the vector field
$$
\Xbold = m_1\left(x_1\,\frac{\partial}{\partial x_2}
                 -x_2\,\frac{\partial}{\partial x_1}\right)
         + \dots
         + m_q\left(x_{2q-1}\,\frac{\partial}{\partial x_{2q}}
                   -x_{2q}\,\frac{\partial}{\partial x_{2q-1}}\right)
$$
 for some positive real numbers~$m_1,\dots,m_q$ where~$q\le (n{+}1)/2$.
(For simplicity, I am taking $x_1,\dots,x_{n+1}$ as an orthogonal linear
set of coordinates on~$\bbR^{n+1}$.) 

By the same analysis as in the previous example, the quotient space~$N$
cannot be defined as a Hausdorff space unless the ratios of the~$m_i$ are
rational, in which case, they can be taken to be integers with greatest
common divisor equal to~1.  Then, the ramification set~$R$ will be larger
than the zero set~$Z$ unless~$m_i=1$ for all~$i$ and it will not be possible 
to define a smooth structure on any extension of~$N$ so that the leaf 
quotient~$\phi:\bbR^{n+1}\setminus R\to N$ extends as a smooth map to
any part of~$R\setminus Z$.  Thus, if the harmonic morphism is to be 
globally defined, one must have~$R=Z$, i.e., all of the~$m_i=1$.  

Now, as has already been noted, the map~$\phi$ cannot be extended smoothly
across~$Z$ unless~$n=3$, so that either~$q=1$ or $q=2$.  With some effort,
it can be shown that the first case leads to a singular quotient that cannot
be smoothed.  However, for~$q=2$, the Killing field
$$
\Xbold =  x_1\,\frac{\partial}{\partial x_2}
         -x_2\,\frac{\partial}{\partial x_1}
         +x_3\,\frac{\partial}{\partial x_4}
         -x_4\,\frac{\partial}{\partial x_3}
$$
does lead to a smooth quotient, for setting
$$
\align
y_1 &= {\ts\frac12}\bigl({x_1}^2+{x_2}^2-{x_3}^2-{x_4}^2\bigr)\\
y_2 &= x_1x_4-x_2x_3\\
y_3 &= x_1x_3+x_2x_4\\
\endalign
$$
yields a map~$\phi = (y_1,y_2,y_3):\bbR^4\to\bbR^3$ that satisfies
$$
\phi^*(h) = |\Xbold|^2\,g - \bigl(\Xbold^\flat\bigr)^2
$$
where~$g$ and $h$ are the standard metrics on~$\bbR^4$ and $\bbR^3$,
respectively.  Thus, this defines a smooth harmonic morphism, globally
defined on~$\bbR^4$.  Note that this is an example of a harmonic morphism
defined by quadratic maps, as studied in~\cite{OW}.  In~\cite{Ba}, it was
shown that this is essentially the only harmonic map from~$\bbR^4$ to~$\bbR^3$
given by quadratic polynomials.
\endexample

\example{Example 4}  {\it Hyperbolic quotients.}
The discussion in this last example will be more cursory than in 
the previous examples since the results are very like the first two
and proved by the same methods.

Again, just as in the Euclidean case, if~$\bigl(M^{n+1},g\bigr)$ is 
the hyperbolic ball, there are two kinds of non-zero Killing vector
 fields.  The first kind have no fixed points, while the second have
a fixed submanifold which is totally geodesic.

If~$\Xbold$ is a Killing field on~$M$ without fixed points, then, just
as in the Euclidean case, one can choose a hypersurface~$N$ in~$M$ that is
transverse to the vector field~$\Xbold$ and so that each integral curve
of~$X$ in~$M$ meets~$N$ in exactly one point.  This~$N$ is diffeomorphic
to~$\bbR^n$, of course, and there will be a unique induced metric on~$N$
so that the leaf projection is a harmonic morphism.

If~$\Xbold$ is a Killing field on~$M$ that has a non-empty fixed point 
set~$Z$, then, again, unless $n=3$, the harmonic morphism defined 
on~$X\setminus Z$ cannot be extended smoothly across the fixed point set.
However, even when~$n=3$, the two cases that arise do not yield a 
smooth quotient, being more like the two cases that arose in the study
of~$S^4$ than the two cases that arose in the study of~$\bbR^3$.  Thus,
none of these give rise to global harmonic morphisms on hyperbolic
$4$@-space.
\endexample

\enddocument